\documentclass[preprint]{elsarticle}
\usepackage{latexsym, epsfig}
\usepackage{amsmath}
\usepackage{amssymb}
\usepackage{amsfonts}
\usepackage{graphicx}
\usepackage{subfig}
\usepackage{adjustbox}
\usepackage{algorithm,algpseudocode}
\usepackage{bm}
\usepackage{textcomp}
\usepackage[thinc]{esdiff}
\usepackage[dvipsnames]{xcolor}
\usepackage{comment}
\usepackage{courier}
\usepackage{tikz}
\usepackage{makecell}
\usepackage{hyperref}

\usepackage{longtable}

\usepackage[toc,page,titletoc]{appendix}
\usepackage[capitalise,nameinlink]{cleveref}

\crefname{supp}{Supplement}{Supplements}

\usepackage{booktabs}
\usepackage{caption}
\usepackage{float}
\usepackage{color,soul}

\newcommand{\tran}{^{\text{T}}} 

\newcommand{\inv}{^{-1}}

\newcommand{\Amat}{\mathbf{A}}
\newcommand{\Gmat}{\mathbf{G}}
\newcommand{\Imat}{\mathbf{I}}
\newcommand{\Mmat}{\mathbf{M}}
\newcommand{\Nmat}{\mathbf{N}}
\newcommand{\Lmat}{\mathbf{L}}
\newcommand{\Dmat}{\mathbf{D}}
\newcommand{\Umat}{\mathbf{U}}

\newcommand{\fvec}{\mathbf{f}}
\newcommand{\xvec}{\mathbf{x}}
\newcommand{\yvec}{\mathbf{y}}
\newcommand{\rvec}{\mathbf{r}}
\newcommand{\evec}{\mathbf{e}}
\newcommand{\uvec}{\mathbf{u}}
\newcommand{\vvec}{\mathbf{v}}
\newcommand{\phivec}{\bm{\phi}}
\newcommand{\zerovec}{\mathbf{0}}

\newcommand{\Gopt}{\mathcal{G}}

\newcommand{\Lopt}{\mathcal{L}}
\newcommand{\Bopt}{\mathcal{B}}
\newcommand{\Oopt}{\mathcal{O}}

\newcommand{\hDON}{h_\text{D}}

\newcommand{\newstep}{^{(k+1)}}
\newcommand{\curstep}{^{(k)}}

\newcommand{\initstep}{^{(0)}}

\newcommand{\Cov}{\text{Cov}}

\newcommand{\Rdomain}{\mathbb{R}}

\algrenewcommand\algorithmicrequire{\textbf{Input:}}
\algrenewcommand\algorithmicensure{\textbf{Output:}}

\appendixpageoff
\appendixtitleoff

\newcommand\revv[1]{#1}
\newcommand\oldtext[1]{}

\renewcommand{\textcolor}[2]{#2}

\journal{Journal}
\begin{document}
\begin{frontmatter}
\title{Blending Neural Operators and Relaxation Methods in PDE Numerical Solvers}

\author[1]{Enrui Zhang \corref{contrib}}
\author[1]{Adar Kahana \corref{contrib}}
\author[1]{Alena Kopani\v{c}\'akov\'a}
\author[2]{Eli Turkel}
\author[3]{Rishikesh Ranade}
\author[3]{Jay Pathak}
\author[1,4]{George Em Karniadakis\corref{cor1}}
\cortext[contrib]{Authors contributed equally}
\cortext[cor1]{Corresponding author E-mail: george\_karniadakis@brown.edu}
\address[1]{Division of Applied Mathematics, Brown University, Providence, RI 02912, USA}
\address[2]{Department of Applied Mathematics, Tel Aviv University, Tel Aviv 69978, Israel}
\address[3]{CTO Office, Ansys Inc, Canonsburg, PA 15317, USA}
\address[4]{School of Engineering, Brown University, Providence, RI 02912, USA}

\begin{abstract}
Neural networks suffer from spectral bias having difficulty in representing the high frequency components of a function while relaxation methods can resolve high frequencies efficiently but stall at moderate to low frequencies. 
We exploit the weaknesses of the two approaches by combining them synergistically to develop a fast numerical solver of partial differential equations (PDEs) at scale. 
Specifically, we propose HINTS, a hybrid, iterative, numerical, and transferable solver by integrating a Deep Operator Network (DeepONet) with standard relaxation methods, leading to parallel efficiency and algorithmic scalability for a wide class of PDEs, not tractable with existing monolithic solvers. 
HINTS balances the convergence behavior across the spectrum of eigenmodes by utilizing the spectral bias of DeepONet, resulting in a uniform convergence rate and hence exceptional performance of the hybrid solver overall. 
Moreover, HINTS applies to \textcolor{red}{large-scale, multidimensional systems}, it is flexible with regards to \textcolor{red}{discretizations, computational domain, and boundary conditions}
\revv{and it can be used to precondition Krylov methods as well.} 
\end{abstract}

\end{frontmatter}

\section{Introduction}
Since the proposal of numerical methods for solving differential equations more than half a century ago, scientists and engineers have been able to significantly expand knowledge and insights that have never been achieved in the analytical era, in all fields of physical sciences and engineering~\cite{chapra2011numerical,mathews1992numerical} such as astrophysics~\cite{bodenheimer2006numerical}, quantum physics~\cite{feit1982solution}, fluid dynamics~\cite{patera1984spectral,kim1987turbulence,cockburn2012discontinuous}, solid mechanics~\cite{hughes2012finite,simo2006computational,hughes2005isogeometric,jing2002numerical}, material science~\cite{rappaz2003numerical}, and electrodynamics~\cite{kong2004scattering}. With an appropriate numerical method such as finite differences~\cite{strikwerda2004finite}, finite elements~\cite{hughes2012finite,bathe2006finite}, or a spectral method~\cite{karniadakis2005spectral}, one may obtain the solution of differential equations that may involve nonlinearity, complex geometry, and/or multiscale phenomena. Thus, one acquires a quantitative knowledge and understanding of the governing mechanism of physical and engineering systems.

Despite the diversity of differential equations in different areas, their numerical solutions, in many cases, are reduced to solving systems of linear equations as the last step. This seemingly simple task is, however, not perfectly resolved yet so far. Existing numerical solvers are frequently convergent and stable only for linear systems satisfying certain conditions (e.g., positive definiteness)~\cite{burden2015numerical}. Computational efficiency becomes an issue when dealing with large-scale systems. Among the iterative solvers, the Jacobi method and the Gauss-Seidel method (see a brief review in Section S1 in Supplementary Information (SI)), as the most classical methods, suffer from divergence for non-symmetric and indefinite systems and slow convergence associated with low-frequency eigenmodes. More advanced iterative methods~\cite{van2003iterative,jinchao_solvers,greenbaum1997iterative,olshanskii2014iterative,saad2003iterative,knoll2004jacobian,gander2022unified,briggs2000multigrid,hackbusch2013multi} and techniques have been proposed to mitigate these two issues. In particular, multigrid methods~\cite{briggs2000multigrid,hackbusch2013multi,bramble2019multigrid,berger2012inexact,hiriyur2012quasi,bramble1990parallel} (see Section S1.2) use a hierarchy of discretizations to efficiently relax eigenmodes of all frequencies, making them one of the state-of-the-art methods in terms of computational efficiency. Despite being a century-old problem, solving linear systems is still an open question, calling for better algorithms that are more efficient and applicable to broader classes of linear systems.

In the recent decade, scientific machine learning (SciML) has developed rapidly for problems in physical sciences and engineering due to its potential advantage in predictive capability and efficiency. In particular, deep learning for operator regression is an active research area of wide interest for system identification and fast inference. Recent works have proposed several frameworks, including the deep operator network (DeepONet)~\cite{lu2021deeponet}, Fourier Neural Operator (FNO)~\cite{fno}, etc.~\cite{kissas2022learning,patel2021physics}. Among these methods, DeepONet was the original method and was proposed based on the universal approximation theorem for operators~\cite{univapproxtheorem}, showing a strong prediction capability in diverse engineering fields~\cite{don_fair,don_bubble,don_mm,don_brittle,don_dissection,don_coupling,don_phase, don_transfer, don_constitutive} (Section S1.3). It provides fast inference for mappings between functions, which can be highly nonlinear and may take hours or even days to simulate using traditional approaches.

In some SciML works, researchers introduce learning algorithms into numerical solvers for differential equations and/or their associated linear systems. Briefly, these works seek to replace and/or improve existing numerical solvers. Physics-informed neural networks~\cite{pinn,pinn_hfm,pinn_fluid,pinn_geometry,pinn_material,pinn_sysbio,pinn_optics} and related works~\cite{pidon,dgm} represent a unique methodology where the differential equation is explicitly encoded into the loss function, hence entirely eliminating the need to employ numerical solvers during both training and testing stages. Among the majority of studies without such explicit encoding, a deep learning model needs to be trained with sufficient data, typically acquired from simulated data, experimental measurements, or a combination of both. Once the training completes, it provides fast and approximate predictions without invoking numerical solvers~\cite{pfaff2020learning,xue2020amortized,kochkov2021machine,long2018pde,kahana2020obstacle,ovadia2021beyond,tompson2017accelerating}. Several other studies focus on improvements on top of numerical solvers, such as: a correction on coarse-grid results~\cite{um2020solver}, a better iterative scheme~\cite{hsieh2019learning,he2019mgnet,chen2022meta}, a better initial guess of the solution~\cite{huang2020int}, a better meshing strategy~\cite{zhang2020meshingnet,kato2018neural}, a better discretization operator~\cite{bar2019learning}, a better prolongator/restrictor~\cite{luz2020learning,greenfeld2019learning} or a preconditioner~\cite{azulay2022multigrid} for multigrid methods.

In this paper, we propose a fundamentally different approach to integrate DeepONet (Fig. S.1) with standard relaxation methods, yielding HINTS -- a hybrid, iterative, numerical, and transferable solver for differential equations.
Our objective is to utilize the merits from both sides, i.e., neural operators and classical solvers, to design a fast, accurate, and widely applicable solver for differential equations. We demonstrate the effectiveness of HINTS and analyze its characteristics by presenting a series of numerical examples, including \textcolor{red}{different choices of integrated standard solvers, differential equations, and spatial geometries.}

\section{Results}
\subsection{HINTS: Integrating DeepONet and Relaxation Solvers}
We consider the following linear differential equation
\begin{subequations}
\label{eqn:diffeqn_general}
\begin{align}
    \label{eqn:diffeqn_general_inside}
    \Lopt_\xvec(u;k)=f, \xvec \in \Omega\\
    \label{eqn:diffeqn_general_bc}
    \Bopt_\xvec(u)=g, \xvec \in \partial\Omega
\end{align}
\end{subequations}
where $\Lopt_\xvec$ is the differential operator, $\Bopt_\xvec$ is the boundary operator, $k=k(\xvec)$ parameterizes $\Lopt_\xvec$, $f=f(\xvec)$ and $g=g(\xvec)$ are the right-hand-side and boundary terms, and $u=u(\xvec)$ is the solution. Assuming that $g$ is a known fixed function, Eq.~\ref{eqn:diffeqn_general} defines a family of differential equations parameterized by $f$ and $k$.

Before employing HINTS, we first need to train a DeepONet offline, see also Fig. S.2.
This DeepONet approximates the solution operator $\Gopt$ defined by
\begin{align}
    \label{eqn:deeponet}
    \Gopt: k, f \mapsto u \text{ s.t. Eq.~\ref{eqn:diffeqn_general} holds},
\end{align}
where we have assumed the uniqueness of the solution $u$. DeepONet receives $k(x)$ and $f(x)$ in the form of discrete evaluations on $n_{\text{D}}+1$ uniform points, respectively. We denote this discretization associated with DeepONet as $\Omega^{\hDON}$. 

\begin{figure}
  \centering
\includegraphics[width=0.75\textwidth]{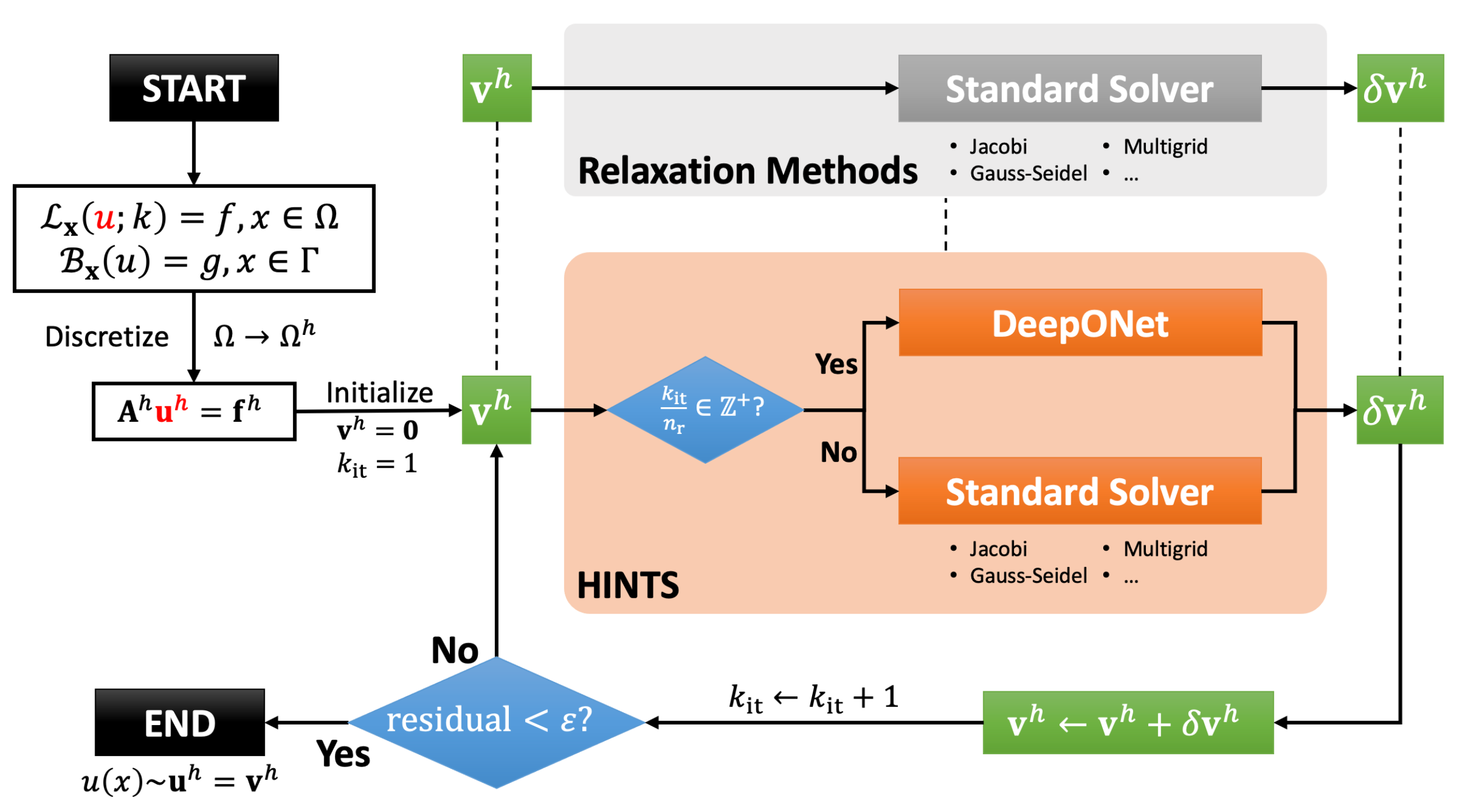}
  \caption{\textbf{Overview of the Hybrid Iterative Numerical Transferable Solver (HINTS).} The goal is to solve $u(\xvec)$ in the differential equation. HINTS starts by discretizing the computational domain $\Omega$ into $\Omega^h$, where $h$ denotes the mesh size.
We then initialize a guess of the solution ($\vvec^h=0$, for example). In the $k_\text{it}$th iteration, the approximate solution $\vvec^h$ is corrected by $\delta\vvec^h$ through either the DeepONet solver or the standard numerical solver. The choice of the solver is determined by whether $n_\text{r}$ divides $k_\text{it}$, where $n_\text{r}$ is a parameter that we choose to optimize the performance. There is a plurality of choices for the numerical solver, such as Jacobi, Gauss-Seidel, and multigrid methods. 
The algorithm proceeds until the residual of the underlying PDE is smaller than a threshold $\varepsilon$.}
  \label{fig:sketch}
\end{figure}

The workflow of a general HINTS is shown in Fig.~\ref{fig:sketch}. After discretizing the differential equation into linear systems, traditional solvers adopt a fixed relaxation method in each iteration. HINTS, on the other hand, creates a second branch of the iterator using a trained DeepONet. During the process of the iterative solution, we alternately adopt the numerical iterator and the DeepONet iterator with ratio $1:(n_{\text{r}}-1)$ (i.e., DeepONet with proportion $1/n_\text{r}$), until the solution converges. The numerical iterator may be chosen from established solvers, such as the Jacobi method, the Gauss-Seidel method, etc.

Among the choices of the numerical iterators, we mainly consider the (damped) Jacobi method, formulating the HINTS-Jacobi (Algorithm S3). For a class of problems defined by Eq.~\ref{eqn:diffeqn_general}, the target is to compute the solution $u(\xvec)$ given $k(\xvec)$ and $f(\xvec)$. The algorithm starts by discretizing the domain $\Omega$ into $\Omega^h$ and assembling the $n$-dimensional linear system $\Amat^h\uvec^h=\fvec^h$. 
In each iteration, given the previous approximate solution $\vvec^h$ ($\approx\uvec^h$), the algorithm seeks to solve $\delta\vvec^h$ in $\Amat^h\delta\vvec^h=\rvec^h$, where the residual $\rvec^h$ is defined as
\begin{align}
    \label{eqn:residual}
    \rvec^h:=\fvec^h-\Amat^h\vvec^h,
\end{align}
and then make a correction to the previous $\vvec^h$ by $\vvec^h\leftarrow \vvec^h+\delta \vvec^h$. To calculate the correction $\delta\vvec^h$, either DeepONet or the Jacobi solver is invoked. For the DeepONet solver, it takes $k(\xvec)$ and $r(\xvec)$ (function-form residual corresponding to the vector-form $\rvec^h$) as inputs, and yields $v(\xvec)=\delta \vvec^h$ as the output. The algorithm proceeds until the residual is sufficiently small. 
\oldtext{By replacing the embedded Jacobi solver with Gauss-Seidel solver, we can similarly formulate the HINTS-GS.}
\revv{We can replace Jacobi with any other relaxation method, giving rise to different variants of the HINTS algorithm. 
For instance, considering Gauss-Seidel relaxation leads to the HINTS-GS method.}
    
Another more advanced option is to consider combining DeepONet with multigrid methods, yielding the HINTS-MG (Algorithm S4). 
Within a multigrid algorithm, we replace a relaxation at each grid level by the HINTS-based relaxation. 
\revv{This is of particular interest for problems, for which the standard relaxation methods do not constitute effective and numerically robust smoothers, such as Helmholtz problems, c.f.~\cite{elman2001multigrid}.}

Notably, the DeepONet solver is flexible with regards to the domain discretization of the numerical problem. A DeepONet trained on discretization $\Omega^{\hDON}$ can be transferred to a numerical problem with another discretization $\Omega^h$, by simply including an interpolation step in the input side of DeepONet. 
This enables a HINTS to be used for solving equations equipped with different discretizations, involving different mesh densities, non-uniform meshes, and/or irregular meshes. 
It is also possible to extend DeepONet to a different domains using the technique for domain adaptation proposed recently in~\cite{don_transfer}.

We employ HINTS for the following two prototypical equations:
\begin{itemize}
    \item Poisson equation:
    \begin{subequations}
    \label{eqn:poisson}
    \begin{align}
        \label{eqn:poisson_PDE}
        \nabla\cdot\Big(k(\xvec)\nabla u(\xvec)\Big)+f(\xvec)&=0, \quad \xvec\in\Omega\subset\mathbb{R}^d\\
        \label{eqn:poisson_BC}
        u(\xvec)&=0, \quad \xvec\in\partial\Omega;
    \end{align}
    \end{subequations}
    \item Helmholtz equation:
    \begin{subequations}
    \label{eqn:helmholtz}
    \begin{align}
        \label{eqn:helmholtz_PDE}
        \nabla^2 u(\xvec)+k^2(\xvec)u(\xvec)&=f(\xvec), \quad \xvec\in\Omega\subset\mathbb{R}^d\\
        \label{eqn:helmholtz_BC}
        u(\xvec)&=0, \quad \xvec\in\partial\Omega.
    \end{align}
    \end{subequations}
\end{itemize}
The Helmholtz equation is indefinite for sufficiently large $k(\xvec)$, making it hard to solve by many classical solvers. In the following subsections, we will examine the performance of HINTS by presenting a series of numerical examples for these two equations, including different spatial dimensions, different discretization techniques, and different geometries.

\subsection{Poisson Equation in One Dimension}
We first consider the prototypical Poisson equation in one dimension, defined in $\Omega=(0,1)$. The goal of HINTS-Jacobi is to solve this equation with arbitrary $k(x)$ and $f(x)$. We first train a DeepONet with paired data $[k(x),f(x)]$ (generated by a Gaussian random field (GRF), as outlined in Section S2.1 and Tab.S.1) and corresponding $u(x)$. 
After training, we employ HINTS-Jacobi to solve for new instances of $k(x)$ and $f(x)$. Here, we employ the same discretization for the DeepONet and the numerical problem ($n_\text{D}=n=30$).

\begin{figure}
  \centering
\includegraphics[width=0.975\textwidth]{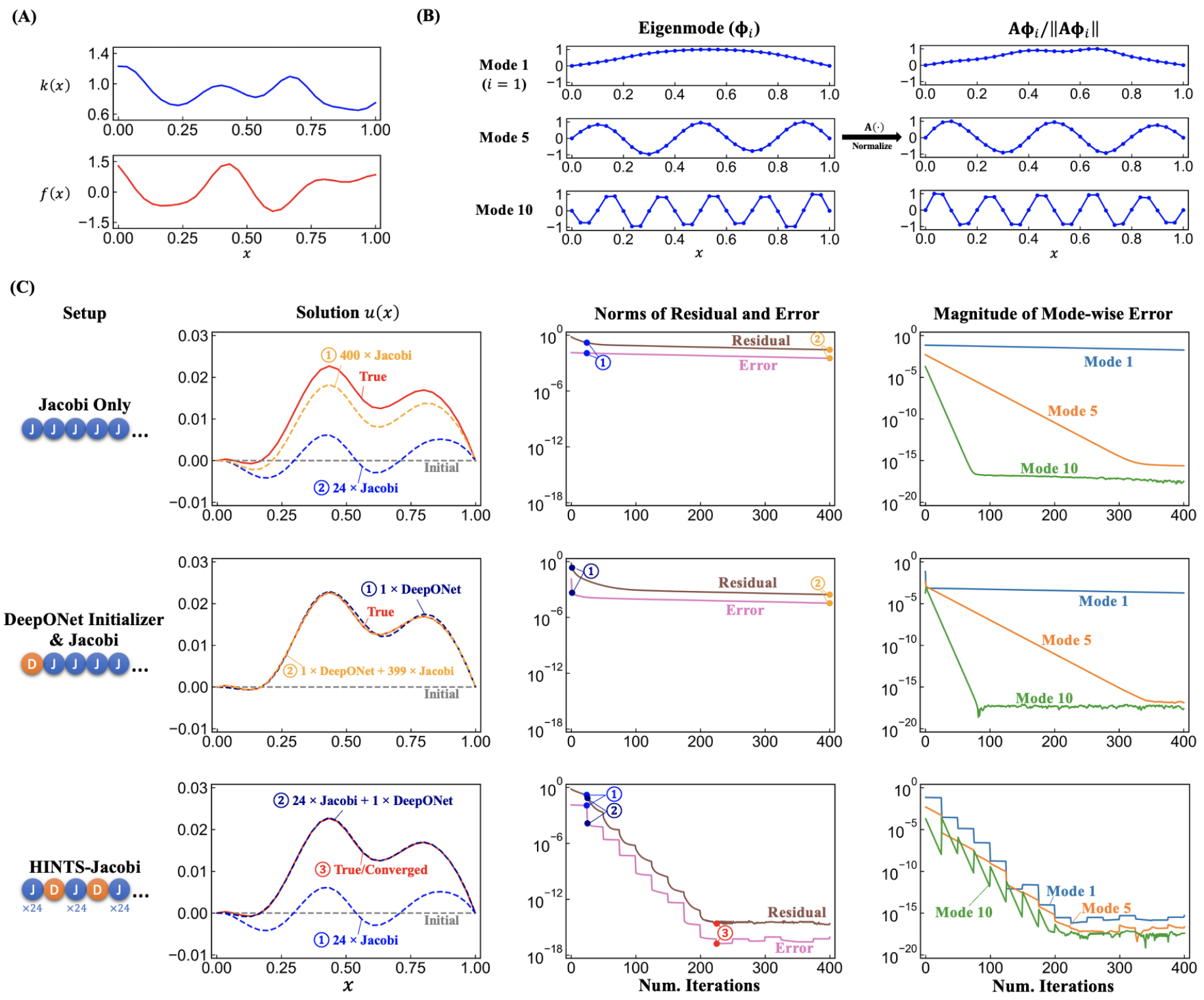}
\caption{\textbf{Results of 1D Poisson Equation.} (A) Profiles of $k(x)$ and $f(x)$. (B) Eigenmodes $\phivec^h_i$ ($\phi_i(x)$; $i=1,5,10$ herein) and corresponding loading vectors $\frac{\Amat^h\phivec^h_i}{||\Amat^h\phivec^h_i||}$. Mode 1 has the lowest spatial frequency, while mode 10 has a relatively high spatial frequency. (C) Numerical results. We consider three setups, each shown in one row: (M1) Jacobi solver only; (M2) Jacobi solver with DeepONet initializer, i.e., one-time usage of DeepONet followed by Jacobi iterations; (M3) HINTS-Jacobi (with a DeepONet-to-Jacobi ratio $1:24$). The second column shows key snapshots of the iterative solution. 
The third column shows the histories of the norms of residual and error of the approximate solution, with the snapshots in the second column marked correspondingly. The fourth column shows the history of the norm of error for eigenmodes 1, 5, and 10.}
  \label{fig:results_1P}
\end{figure}

We show the performance of HINTS-Jacobi for a representative choice of $(k(x),f(x))$ in Fig.~\ref{fig:results_1P}. 
Results for other $(k(x),f(x))$ and for HINTS-GS are presented in Section S3.2. 
Fig.~\ref{fig:results_1P}(A) shows the profiles of $k(x)$ and $f(x)$ considered. In Fig.~\ref{fig:results_1P}(B), we show the profiles of typical eigenmodes $\phivec^h_i$ ($\phi_i(x)$; $i\in\{1,5,10\}$) of the solution $u(x)$ and corresponding loading vectors ($\frac{\Amat^h\phivec^h_i}{||\Amat^h\phivec^h_i||}$), where mode 1 refers to the lowest-frequency and mode 10 is a  relatively high frequency. In Fig.~\ref{fig:results_1P}(C), we present the numerical results for three setups: (M1) Jacobi solver; (M2) Jacobi solver initialized with a single DeepONet inference; (M3) HINTS-Jacobi ($1/25$).
We show key snapshots of the approximate solution and the norms of residual and error of the approximate solution, where the reference solution is obtained using a direct solver.
After 400 iterations, the approximate solution of M1 is still visually different from the true solution $u(x)$, with the absolute error being $\Oopt(10^{-2})$. 
M2 provides a better estimate, with an absolute error $\Oopt(10^{-4})$ that is substantially smaller than M1. However, the contribution of the one-time DeepONet initialization in M2 is marginal -- it does not change the convergence rate, and the solution is still far from the convergence. 
On the contrary, M3, i.e., HINTS-Jacobi, provides a solution that rapidly converges to machine zero within roughly 200 iterations. 

In the last column of Fig.~\ref{fig:results_1P}(C), we display the history of the eigenmode-wise error for modes 1, 5, and 10. For all the setups, the convergence rate of the embedded Jacobi iterations is uniform for each mode, matching the theoretical analysis in Section S1.1. 
For Setups M1 and M2, the slow convergence of mode 1 in the Jacobi solver is the bottleneck of the overall convergence rate. The embedded DeepONet iterations, on the other hand, make a substantial difference in M3. DeepONet iterations significantly reduce mode-1 error, hence breaking the bottleneck effect in restricting the overall convergence rate. While DeepONet iterations have limited improvements for mode 5 and even negative effects on mode 10, these two modes are well resolved by the Jacobi solver. 
With such an observation, we find that the mixture of the two solvers functions in a synergetic way: \emph{the Jacobi solver is efficient for high-frequency modes but not for low-frequency modes; the DeepONet solver provides fast approximate solutions for low frequencies but may pollute the high frequencies. A proper combination of these two methods synchronizes the convergence paths and enables fast and uniform convergence across all eigenmodes.}

\subsection{Solving Equations Involving Indefiniteness}
\begin{figure}
  \centering
\includegraphics[width=0.975\textwidth]{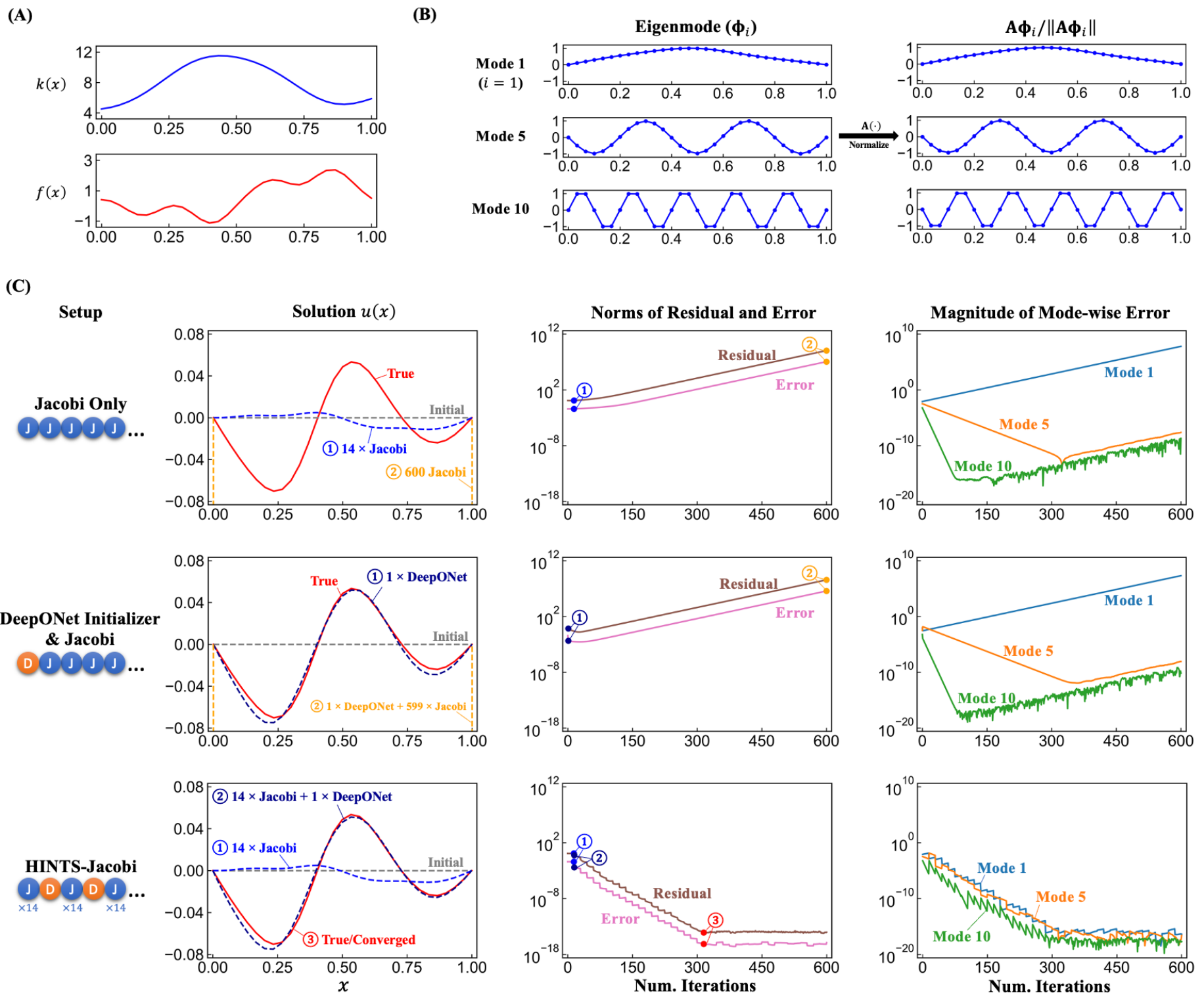}
  \caption{\textbf{Results of 1D Helmholtz Equation.} (A) Profiles of $k(x)$ and $f(x)$. (B) Eigenmodes $\phivec^h_i$ ($\phi_i(x)$; $i=1,5,10$ herein) and corresponding loading vectors $\frac{\Amat^h\phivec^h_i}{||\Amat^h\phivec^h_i||}$. 
  Mode 1 has the lowest spatial frequency, while mode 10 has a relatively high spatial frequency. (C) Numerical results. We consider three setups, each shown in one row: (M1) Jacobi solver only; (M2) Jacobi solver with DeepONet initializer, i.e., one-time usage of DeepONet followed by Jacobi iterations; (M3) HINTS-Jacobi (with a DeepONet-to-Jacobi ratio $1:14$). The second column shows key snapshots of the iterative solution. The third column shows the histories of the norms of residual and error of the approximate solution, with the snapshots in the second column marked correspondingly. The fourth column shows the history of the norm of error in eigenmodes 1, 5, and 10.}
  \label{fig:results_1H}
\end{figure}

In the preceding example of the Poisson equation, the linear system is positive definite, which can be solved by majority of  classical numerical solvers. Here, we consider the 1D Helmholtz equation (Eq.~\ref{eqn:helmholtz}) in $\Omega=(0,1)$, where the indefiniteness causes a divergence of many classical solvers. After the offline training of a DeepONet, we use the HINTS-Jacobi to solve for new instances $k(x)$ and $f(x)$ with the same discretization factor $n=n_\text{D}=30$. 
We display the results in Fig.~\ref{fig:results_1H}.
As analyzed in Section S1.1, the first two setups diverge due to the divergence of low-frequency modes for indefinite systems. 
Setup M3, i.e., the HINTS-Jacobi ($1/15$), is capable of converging after roughly 300 iterations. Here, HINTS-Jacobi functions in a similar but not identical way as for Poisson equation -- the DeepONet solver decreases the error of low-frequency modes, counteracting the Jacobi solver that slowly increases the error, while leaving high-frequency modes to be handled by the Jacobi solver. The HINTS-Jacobi resolves the issue of divergence for indefinite systems, \revv{caused by amplification of the low-frequency error components}, which is prevalent for the classical stationary methods.

\subsection{Applicability to higher dimensions and irregular geometries}
\begin{figure}
  \centering
\includegraphics[width=0.9\textwidth]{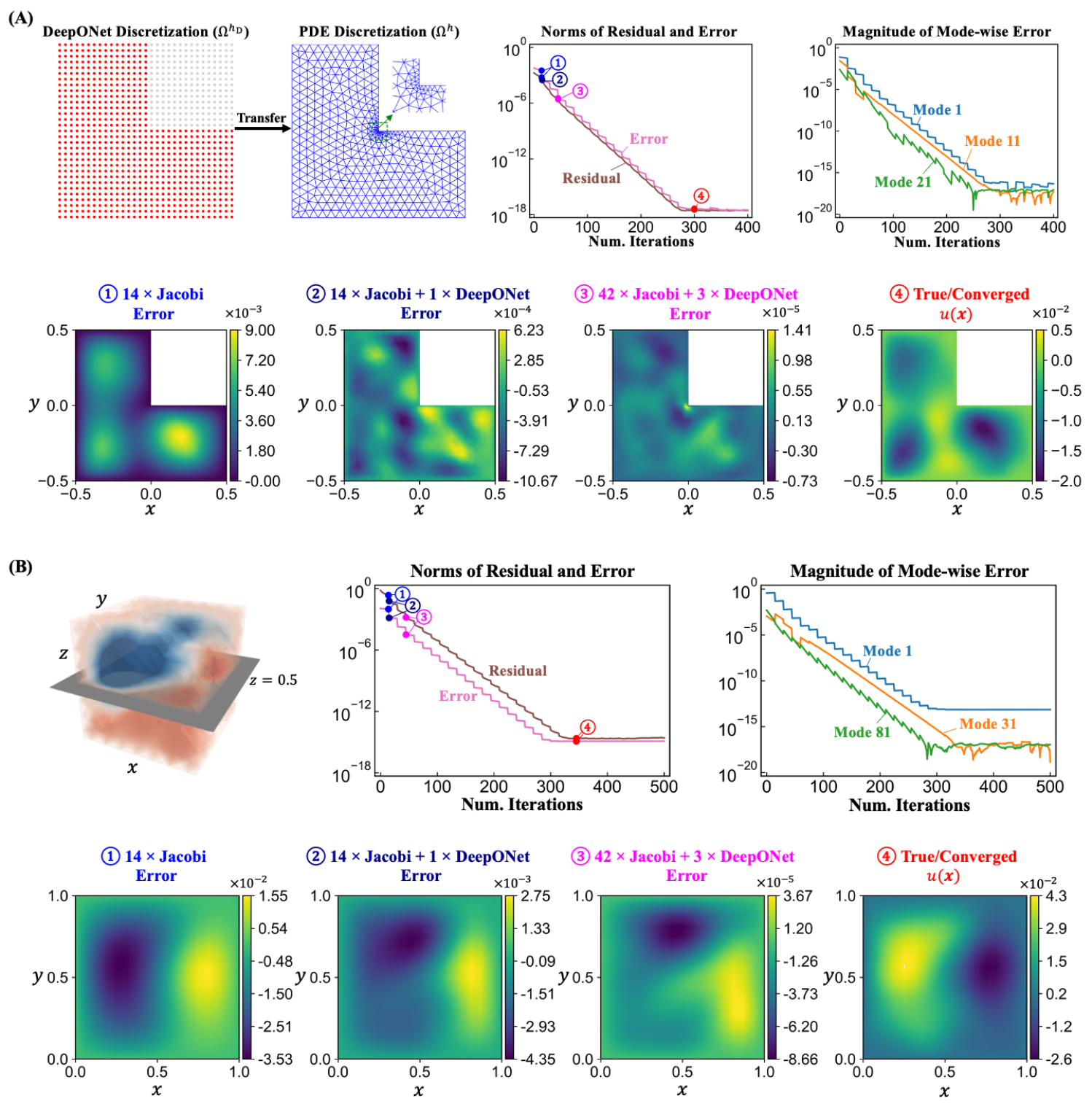}
\caption{\textbf{Results of HINTS-Jacobi for 2D Poisson Equation and 3D Helmholtz Equation.} The DeepONet-to-Jacobi ratio is $1:14$ in both cases. (A) 2D Poisson equation. In the first row, the first and the second columns shows the discretization of DeepONet and the numerical problem, respectively; the third and the fourth columns show the histories of norms of errors and residuals and mode-wise errors, respectively. In the second row, we display the profiles of errors of some key snapshots together with the true/converged solution. (B) 3D Helmholtz equation. In the first row, the first column illustrates the 3D domain where we define the problem. Other displays are similar to Poisson equation in (A). The snapshots show the errors and true solution on cross section $z=0.5$, marked in the first panel.}
  \label{fig:results_HD}
\end{figure}

The effectiveness of HINTS is not limited to one-dimensional problems equipped with uniform grids. 
Due to the discretization transferability, HINTS is capable of adapting to an irregular computational domain (e.g., triangular $\Omega$) and/or irregular discretizations (e.g., $\Omega^h$ with triangular elements). 
Here, we consider the following examples:
\begin{itemize}
    \item 2D Poisson equation defined in an L-shaped domain: $\Omega:=(0,1)^2\backslash[0.5,1)^2$;
    \item 3D Helmholtz equation defined in an unit square: $\Omega:=(0,1)^3$.
\end{itemize}
We illustrate the geometries in the first column of the first rows in Fig.~\ref{fig:results_HD}(A) (2D Poisson) and Fig.~\ref{fig:results_HD}(B) (3D Helmholtz). In particular, the case of 2D Poisson equation uses finite element discretization with linear triangular meshes ($\Omega^h$), while we still use a DeepONet trained on a uniform mesh ($\Omega^{\hDON}$). Both cases adopt the DeepONet proportion $1/15$. 
The remaining results in the first rows of Figs.~\ref{fig:results_HD}A-B, similar to 1D cases, show the histories of the norms of the residual, error, and the mode-wise error. The second rows of Figs.~\ref{fig:results_HD}A-B show the error of three key snapshots during the iterations together with the true/converged solution. For the case of 3D Helmholtz equation, we display the cross section at $z=0.5$. In both cases, DeepONet iterations significantly reduce the errors or low-frequency modes, hence accelerating the convergence of the solution. 
Figures S.7(A-C) present additional results of a 2D Helmholtz equation defined in $\Omega:=(0,1)^2$ and 3D Helmholtz equation solved with HINTS-GS.
Moreover, Fig. S.8 provides an illustration of eigenmodes for problems with varying geometries.

\subsection{Integrating HINTS with Multiscale Methods}
\begin{figure}
  \centering
 \includegraphics[width=\textwidth]{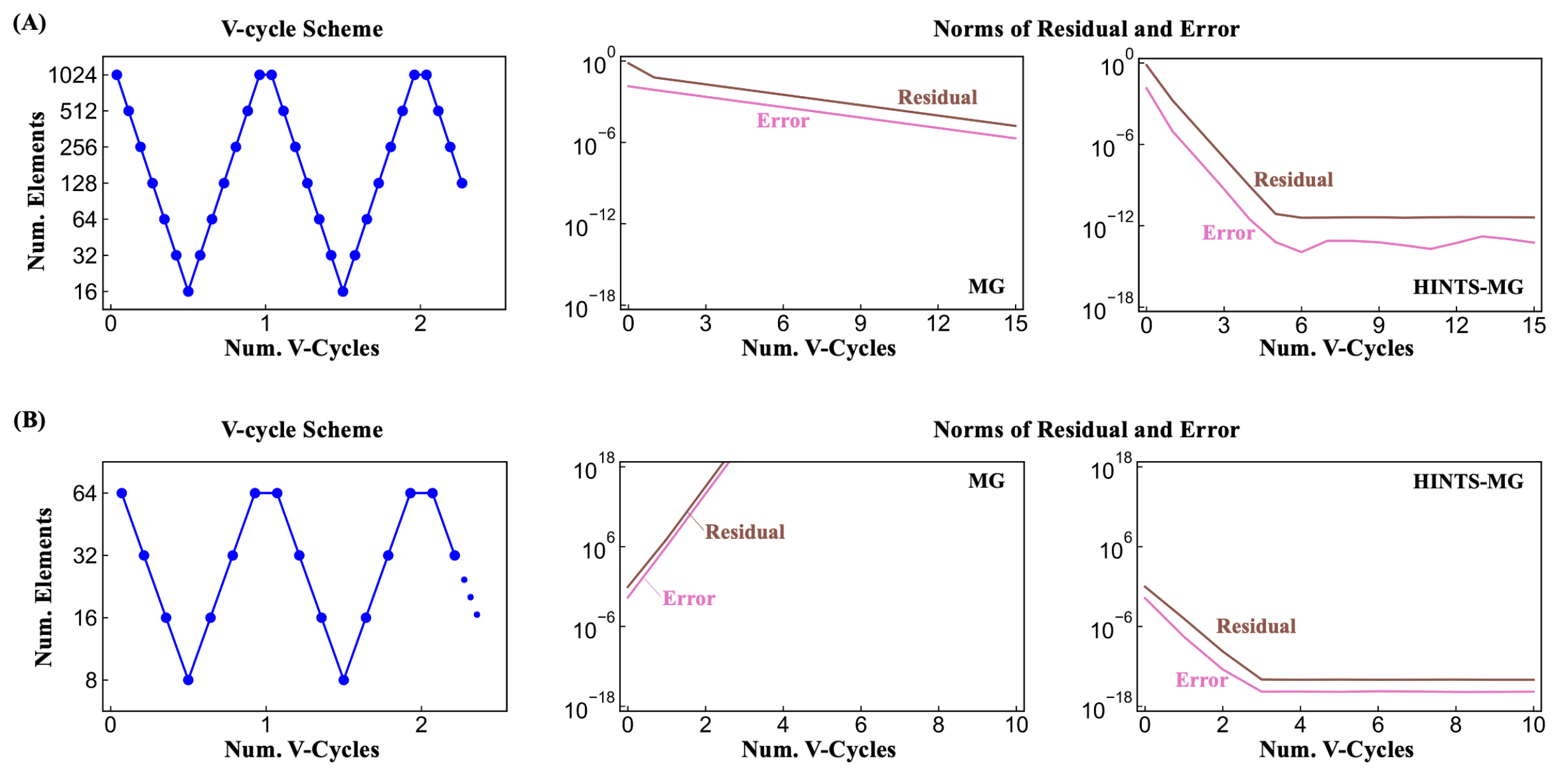}
\caption{\textbf{Results of HINTS-MG for 1D Poisson Equation and 2D Helmholtz Equation.} (A) 1D Poisson Equation. The first column shows the setup of V-cycle in the multigrid algorithm, where each dot represents 10 iterations. The second and the third column compare the histories of norms of errors and residuals for the two setups of relaxation: (M1) with Gauss-Seidel relaxation only; (M2) with hybrid DeepONet-GS relaxation ($1:9$ in each grid level). (B) 2D Helmholtz equation. Display is similar to (A). The DeepONet-GS ratio is also $1:9$.}
  \label{fig:results_MG}
\end{figure}

\revv{The DeepONet augmented iterative methods can naturally be employed as subspace solvers for large-scale multiscale methods~\cite{xu1992iterative}, e.g., multigrid or domain-decomposition methods. This is of particular interest for problems for which efficient, robust, and numerically stable subspace solvers are not straightforward to construct using standard numerical approaches.}
We demonstrate such capability of HINTS by devising the hybrid solver that integrates DeepONet with multigrid methods. We apply HINTS-MG to solve the 1D Poisson equation and the 2D Helmholtz equation, with the results shown in Fig.~\ref{fig:results_MG}(A) and Fig.~\ref{fig:results_MG}(B), respectively. 
For the 1D Poisson equation, we consider the same $k(x)$ and $f(x)$ as shown in Fig.~\ref{fig:results_1P}(A). We discretize the domain into $n=1024$ uniform intervals, and adopt the V-cycle scheme with $7$ levels (see Fig.~\ref{fig:results_MG}(A)). We consider two setups for relaxation: (M1) Gauss-Seidel relaxation only; (M2) HINTS-MG with DeepONet-GS relaxation (DeepONet proportion $1/10$). We present the histories of the norms of residual and error (second column) and the mode-wise error (third column) against the number of V-cycles in Fig.~\ref{fig:results_MG}(A). By introducing DeepONet, the number of V-cycles needed for achieving convergence is significantly decreased to roughly 5. For the 2D Helmholtz equation, we show the results for $n=64^2$ and four grid levels in Fig.~\ref{fig:results_MG}(B). While the multigrid solver \revv{with standard Jacobi relaxation} diverges, \revv{due to amplification of low-frequency components of the error on all levels},  HINTS-MG converges within approximately three V-cycles. \textcolor{red}{Notably, due to the capability of HINTS in handling different discretizations at the training and testing stages (i.e., $\Omega^{\hDON}\neq\Omega^h$), one single DeepONet may be used to handle all the grid levels.} 

\subsubsection{Generalization capability of the HINTS}
To investigate the applicability of HINTS for large-scale engineering problems with practical importance, we need to address its generalization capabilities.
Here we summarize the main aspects of the HINTS generalization, and include the corresponding results in Section S3.3.
\begin{itemize}
    \item \textbf{Generalization to different discretizations.} We have shown that the HINTS can be transferred to solve equations where the domain discretization is different from that in the training stage, i.e., $\Omega^{\hDON}\neq\Omega^h$. 
See also Fig. S.5 for extended results concerning mismatching discretizations.
We comment that the generalization capability across different discretizations of HINTS is significant for the adaptability and practical utility of HINTS, extending beyond the results shown in Fig.~\ref{fig:results_HD}(A). 

	\item \textbf{Generalization to out-of-distribution input functions.} We tested the performance of HINTS using unseen function instances drawn from the same GRF distribution as the training data. In Section S3.3, we report the performance of HINTS with out-of-distribution functions, where the GRF has different correlation lengths. Results indicate that the HINTS can well generalize to higher correlation lengths; for smaller correlation lengths, the convergence is slightly slower than the former case, but still significantly faster than without employing DeepONet. We attribute this to the accurate prediction of the DeepONet for low-frequency modes, as analyzed in Section S3.1 and Fig. S.3.

	\item \textbf{Generalization to different computational domains and boundary conditions.} The HINTS applies to irregular computational domains and boundary conditions other than zero Dirichlet, as shown in Fig.~\ref{fig:results_HD} and Fig. S.6, respectively. Moreover, we find that a HINTS trained with a fixed boundary condition and within one computational domain may be generalized to another similar but not identical domains, with an altered boundary condition, without re-training the DeepONet. 
\end{itemize}

\subsubsection{Large Systems and Preconditioning of Krylov Methods}
\label{sec:krylov}
Further, we investigate the capabilities of HINTS for solving large-scale problems $(n > \mathcal{O}(10^6))$.
Notably, due to the generalization capability of discretizations, we still employ a lightweight HINTS, i.e.,~with a DeepONet trained using coarse meshes.
Fig. S.4(B) shows that the HINTS solvers outperform the direct solver for large-scale systems arising from the discretization of the 3D Poisson equation in terms of the wall-clock time and scaling behavior.

Next, we investigate the potential of HINTS solvers to enhance the convergence and the scalability of Krylov methods through preconditioning.  
To this aim, we integrate the HINTS into the software package PETSc~\cite{balay2019petsc}, see~\cite{kopanivcakova2024deeponet, hints_precond_code} for details. 
This integration demonstrates the ability of the HINTS to integrate into existent, well-established software packages. 
Moreover, it allows us to ``blend" HINTS with several readily available large-scale, state-of-the-art solution strategies, giving rise to a wide range of HINTS-based solvers/preconditioners, e.g., HINTS-Jacobi, HINTS-ILU, or HINTS-MG.

\begin{figure}
\includegraphics[width=\textwidth]{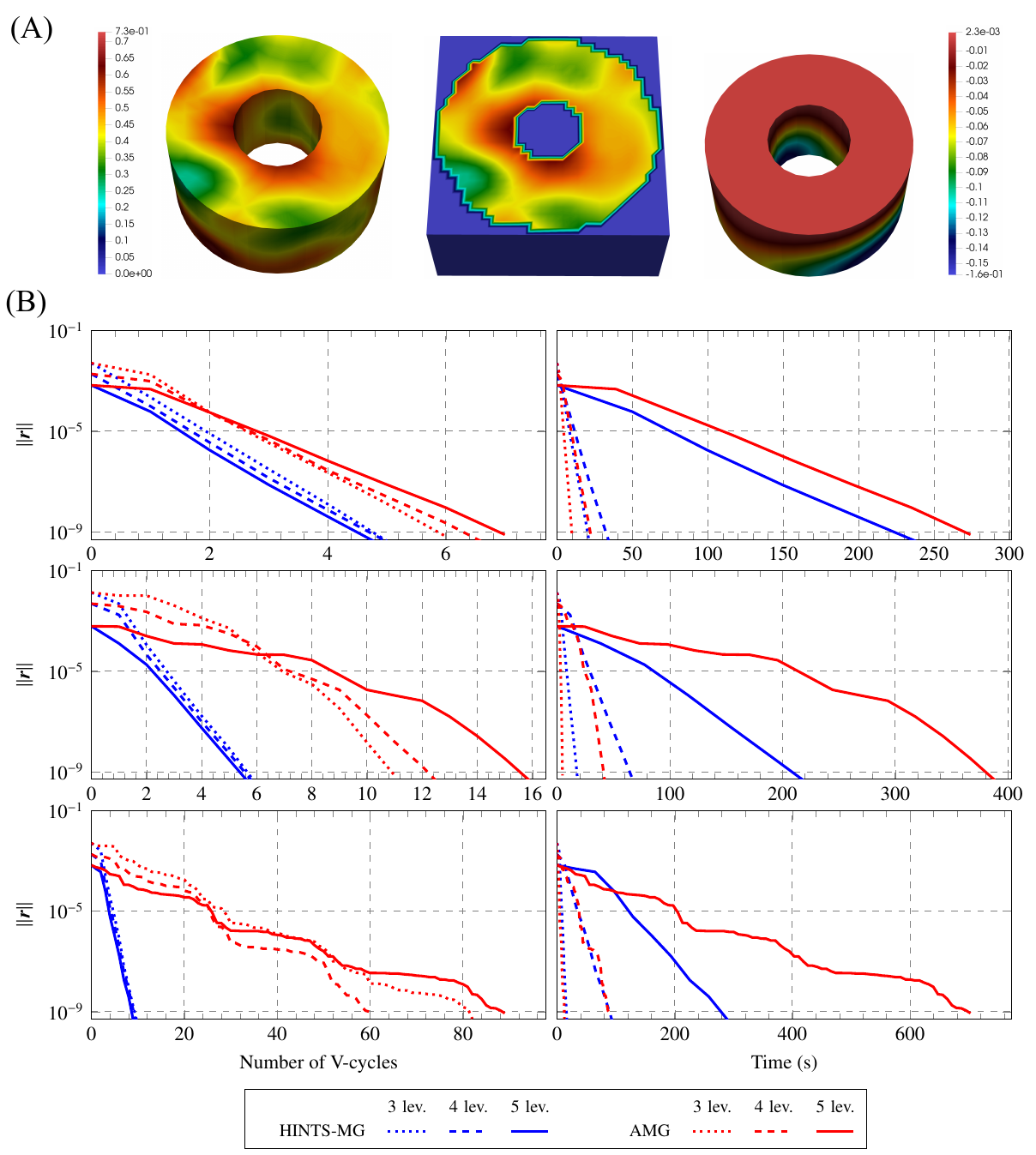}
\caption{\textbf{Results of large-scale Helmholtz 3D - annular cylinder example.}
(A) From left to right: An example of generated right-hand side features, embedding of features into the uniform mesh considered by the DeepOnet and the associated solution. 
(B) Convergence behavior and computational time of F-GMRES method preconditioned with HINTS-MG and AMG in default settings with $k \approx 1, 3$, and $6$ (from top to bottom).
The MGs are configured with three/four/five levels ($74,712$/$570,960$/$4,462,368$ dofs).}
\label{fig:helm3D_cylinder_example}
\end{figure}

We solve the Helmholtz problem~(Eq.~\eqref{eqn:helmholtz_PDE}) using the unstructured geometry of an annular cylinder, see Fig.~\ref{fig:helm3D_cylinder_example}(A) using F-GMRES~\cite{saad1993flexible}. 
Fig.~\ref{fig:helm3D_cylinder_example}(B) on left demonstrates the performance of HINTS-MG preconditioner, which is configured with a varying number of levels and HINTS-Jacobi relaxation.
As we can observe, the F-GMRES with the HINTS-MG preconditioner converges with a fixed number of iterations, independently of the refinement level. 
In contrast, utilizing the algebraic multigrid (AMG) preconditioner, provided by the HYPRE~\cite{falgout2002hypre} and used in its default configuration, results in a slight increase in the number of F-GMRES iterations with increased refinement levels. 
Moreover, as the value of $k$ increases, the convergence of the AMG deteriorates.
Here, we note that several specialized AMG methods for achieving algorithmic scalability for the Helmholtz problem have been proposed in the literature, e.g.,~\cite{olson2010smoothed}. 
However, our objective is not to employ the most optimized solution strategy for a given problem but rather to utilize a state-of-the-art black-box approach. 
This aligns with the main concept behind the proposed HINTS-based solvers/preconditioners, which use DeepONet to augment existent standard solution strategies with the aim of enhancing their convergence, robustness, and algorithmic scalability.

\revv{Fig.~\ref{fig:helm3D_cylinder_example}(B) on right compares the computational time required while using both preconditioners. 
As we can see, the HINTS-MG is more costly for small-scale problems ($n < 0.5 \cdot 10^6$). 
However, as soon as the problem size increases and the problem becomes more challenging ($k$ increases), the benefits of the accelerated convergence are clearly reflected in the obtained timings, and the HINTS-MG preconditioner outperforms the HYPRE-MG by a significant margin.
The details regarding the shortcomings of our current implementation of the HINTS preconditoner in PETSC is analyzed in Section S3.3}

\section{Discussion}
We have demonstrated the capability of HINTS in solving linear differential equations. HINTS is fast, accurate, and widely applicable to diverse classes of differential equations. The method is agnostic in terms of the differential equations, computational domains (shape and dimension), and discretization.

We conducted a systematic study on the performance of HINTS-Jacobi, where we examined why and how this hybrid solver functions. In summary, the embedded DeepONet solver provides fast, approximate solutions for low-frequency eigenmodes in the solution, while the embedded \oldtext{Jacobi solver} \revv{relaxation method} handles the remaining higher-frequency eigenmodes, guaranteeing the convergence of the solver. Such a combination, instead of substitution, incorporates advantages from both ends, improving the performance beyond individual solvers. A hybrid methodology, instead of starting from scratch, retains the advantages from numerical solvers that have been studied for centuries \textcolor{red}{as well as from DeepONet as a powerful operator approximator. 
In Section S3.5, we show inferior results where the DeepONet is replaced by the k-nearest neighbor with Gaussian kernel}. While the DeepONet is trained with a finite dataset, the setup of HINTS enables the generalization into infinite test cases with the preservation of convergence to machine zero.

When training a deep neural network, one frequently observes a phenomenon called spectral bias~\cite{rahaman2019spectral}. This refers to the fact that neural networks tend to learn the lower part of the spectrum first. 
The vast majority of applications of neural networks suffer from such a spectral bias, which inhibits the networks from learning complex mappings. 
On the contrary, HINTS exploits such a bias to tackle low-frequency modes that are otherwise difficult for classical solvers. \oldtext{which, essentially, are biased towards high frequencies.By replacing the classical solver with the DeepONet solver for a limited proportion of iterations, HINTS balances the convergence rates across the spectrum of eigenmodes, significantly improving the computational efficiency and simplifying tremendously implementation and software.}
\revv{Indeed, leveraging spectral bias holds considerable potential for designing robust and scalable large-scale solvers. 
For example, in the context of AMG and domain-decomposition methods, constructing coarse space often requires knowledge of (near-)null-space vectors, a task that may not be well understood in practice for certain classes of problems, such as those with strong heterogeneities. As recently demonstrated in~\cite{kopanivcakova2024deeponet}, utilizing spectral bias, DeepONet can be used to construct the (near-)null-space vectors at low computational cost, in turn facilitating the construction of computationally efficient, scalable, and robust coarse spaces.}

HINTS is promising in practical applications in many fields of real-world science and engineering. One may build HINTS based on existing numerical solvers, by adding a simple conditional statement to integrate DeepONet in a plug-and-play style. \textcolor{red}{Upon the completion of a one-time, offline training procedure, the generalizability of HINTS enables its application for accelerating the solution process of a broad array of equations, with different parameters (e.g., source term, boundary conditions), discretizations, and even computational domains. 
Such a wide applicability effectively renders the training cost nearly insignificant.}
\oldtext{Furthermore, since DeepONet fits into the regular solver routine, HINTS is fully parallelizable. 
In particular, HINTS-Jacobi is embarrassingly parallel, as the update of nodal values requires communication with the neighbor nodes from the preceding iteration, rather than from the current one.
Furthermore, the parallelizability of HINTS endows it with unique advantages in multidimensional, large-scale systems with practical interest, making HINTS \textcolor{red}{a fast} parallel code for problems at scale.}
Furthermore, since DeepONet fits into the regular solver routine, the parallelization of HINTS depends solely on the parallelization capability of the employed iterative method. 
For instance, HINTS-Jacobi is inherently parallel, as the Jacobi method is embarrassingly parallel, and the DeepONet contribution can be evaluated at each computing device for the required nodal point.
This makes HINTS especially well-suited for solving multidimensional, large-scale systems of practical interest.

HINTS represents a general class of frameworks, where one integrates operator regressors with standard numerical solvers. While we have systematically studied several instances, future works are still needed for further investigation of this hybrid methodology. First, one may investigate other solvers to be integrated on both sides. \textcolor{red}{Second, variability regarding geometric domain, boundary conditions, etc., needed further investigation for different problems of interest}. Third, a similar methodology may be extended to nonlinear problems, which are of more practical importance and are more difficult to solve. Finally, the offline cost of training DeepONet can be greatly reduced by transfer learning as we tackle different problems on diverse geometric domains or set of parameters.

\section{Methods}
\label{sec:methods}
This section provides a brief overview of the employed methods. 
Additional technical details are summarized in Section S1.

\subsection{Numerical Solvers}
In the study related to the Jacobi solver, we use the damped Jacobi method, with a damping parameter $\omega=2/3$ for one-dimensional cases, $\omega=4/5$ for two-dimensional cases, $\omega=6/7$ for three-dimensional cases. These choices are shown to be optimal for Poisson equation~\cite{briggs2000multigrid}. 
In the study related to multigrid methods, we consider V-cycles with Gauss-Seidel relaxation \revv{for all tests, except for the Helmholtz-cylinder example, where Jacobi relaxation is employed.}
\oldtext{We conduct Gauss-Seidel relaxation in each grid level (including both the finest and coarsest levels). }
We implemented the full-weighting version of the restrictor and the prolongator.

\subsection{DeepONet}
We use the standard architecture of DeepONet~\cite{lu2021deeponet} in this study. Here, $k(\xvec)$ (similarly for $f(\xvec)$) is fed into the branch network in the discrete form $[k(\xvec_0),k(\xvec_1),...,k(\xvec_{n_{\text{D}}})]\tran$, where $\{\xvec_i\}_{i=0}^{n_{\text{D}}}$ are uniform grid points in $\Omega$. Then, for an arbitrary $\xvec\in\Omega$, which is received at the trunk network, DeepONet predicts $u(\xvec)$ that approximately satisfies Eq.~\eqref{eqn:diffeqn_general}. The branch network is a dense neural network for problems in 1D and a convolutional neural network for problems in multidimensions. We enforce the Dirichlet boundary conditions by postprocessing DeepONet outputs for all cases except the 2D Poisson equation in an L-shaped domain. Specifically, for problems defined in $\Omega=(0,1)^d$ ($d\in\{1,2,3\}$), the original output of the DeepONet is multiplied by $x(x-1)$, $xy(x-1)(y-1)$, $xyz(x-1)(y-1)(z-1)$ for $d=1,2,3$, respectively. To ensure convergence to machine zero, we utilize the scaling of $f(x)$ in the DeepONet. Before feeding into the DeepONet, $f(x)$ is normalized. Then, the norm is multiplied by the output of the DeepONet.

\subsection{HINTS}
The pseudocode related to HINTS is given in Section S1. 
Note that DeepONet predicts the correction of the solution corresponding to the current residual (see Eq.~\ref{eqn:residual}). To solve a single differential equation with a specified $k(\xvec)$ and $f(\xvec)$, DeepONet iterations take the same $k(\xvec)$ as the first input in Eq.~\ref{eqn:deeponet} but different right-hand-sides terms $r(\xvec)$ corresponding to the current residual in the iterative process as the second input. Before invoking the DeepONet solver, we need the conversion of the matrix- (vector-) form residual $\rvec^h$ into the function-form residual $r(\xvec)$. In addition, for numerical examples involving different discretizations of DeepONet and the numerical problem (i.e., $\Omega^h\neq\Omega^{\hDON}$), we need to interpolate $\Omega^h$ into $\Omega^{\hDON}$ before feeding the data into DeepONet. No additional steps are needed on the output side, since DeepONet is capable of evaluating the output function for $\Omega^h$.

\setcounter{equation}{0}
\renewcommand\theequation{S\arabic{equation}}

\setcounter{figure}{0}
\renewcommand\thefigure{S\arabic{figure}}

\setcounter{table}{0}
\renewcommand\thetable{S\arabic{table}}

\setcounter{algorithm}{0}
\renewcommand\thealgorithm{S\arabic{algorithm}}

\renewcommand{\appendixname}{}
\renewcommand{\thesubsection}{S\arabic{subsection}}

\appendix
\setcounter{section}{18}


\section{Supplementary Information}
\subsection{Supplementary Methods}
\label{SI:numerical_methods}

\subsubsection{Jacobi and Gauss-Seidel Methods}
\label{SI_sub:JGS}

We briefly review Jacobi and Gauss-Seidel methods. Through diverse discretization methods (such as finite difference and finite element), solving differential equations can be converted to solving linear equations. We seek $\uvec\in\mathbb{R}^{n}$ in the linear system $\Amat\uvec=\fvec$ ($\Amat\in\mathbb{R}^{n\times n}$ is invertible, $\fvec\in\mathbb{R}^{n}$). Splitting $\Amat=\Mmat+\Nmat$ where $\Mmat$ is invertible, we get $\uvec=-\Mmat\inv\Nmat\uvec+\Mmat\inv\fvec$. This yields the iteration
\begin{equation}
    \label{eqn:iterative_general}
    \uvec\newstep=-\Mmat\inv\Nmat\uvec\curstep+\Mmat\inv\fvec.
\end{equation}
The Jacobi method and the Gauss-Seidel method define $\Mmat$ and $\Nmat$ as
\begin{align}
    \label{eqn:iterative_MN_J}
    \text{Jacobi: } &\Mmat=\Dmat, \Nmat=\Lmat+\Umat,\\
    \label{eqn:iterative_MN_GS}
    \text{Gauss-Seidel: } &\Mmat=\Dmat+\Lmat, \Nmat=\Umat,
\end{align}
where $\Lmat$, $\Dmat$, $\Umat$ are the strictly lower triangular, diagonal, and strictly upper triangular components of $\Amat$, respectively, so that $\Amat=\Lmat+\Dmat+\Umat$. For the Jacobi method, the weighted (damped) version of Eq.~\ref{eqn:iterative_general} is often adopted for better convergence properties, which is defined by
\begin{subequations}
\label{eqn:iterative_weighted}
\begin{align}
    \label{eqn:iterative_weighted_general}
    \uvec\newstep&=(1-\omega)\uvec\curstep-\omega\Mmat\inv\Nmat\uvec\curstep+\omega\Mmat\inv\fvec\\
    \label{eqn:iterative_weighted_J}
    &=(1-\omega)\uvec\curstep-\omega\Dmat\inv(\Lmat+\Umat)\uvec\curstep+\omega\Dmat\inv\fvec,
\end{align}
\end{subequations}
where $\omega$ is the damping parameter. The algorithm is shown in Alg.~\ref{alg:JGS}. During the iterations, for an approximate solution $\uvec\curstep$, we define the residual and the error as
\begin{subequations}
\label{eqn:defs_res_error}
\begin{align}
    \label{eqn:def_residual}
    \rvec\curstep&=\fvec-\Amat\uvec\curstep,\\
    \label{eqn:def_error}
    \evec\curstep&=\uvec\curstep-\uvec,
\end{align}
\end{subequations}
respectively, where $\uvec$ is the true solution and $\rvec\curstep$ and $\evec\curstep$ satisfy $\Amat\evec\curstep=\rvec\curstep$. The evolution of error is given by
\begin{subequations}
\begin{align}
    \evec\newstep &= [(1-\omega)\Imat+\omega(\Imat-\Mmat\inv\Amat)]\evec\curstep\\
    \label{eqn:error_amplification}
    &:= \Gmat_\omega\evec\curstep,
\end{align}
\end{subequations}
where $\Gmat_\omega$ is the amplification matrix, which characterizes how the error vector evolves through iterations.

\begin{algorithm}[h!]
\caption{The (Damped) Jacobi Solver and the Gauss-Seidel Solver}
\label{alg:JGS}
\begin{algorithmic}

\Function{Numerical\_Solver}{$\Amat^h,\fvec^h$}
\State $\vvec^h\gets\zerovec^h$ \Comment{initial guess of the solution}
\State $k_\text{it}\gets 1$
\While{$k_\text{it}\leq n_\text{it}$ and not converged}
    \State $\rvec^h\gets\fvec^h-\Amat^h\vvec^h$
    \If{using damped Jacobi}
        \State $\vvec^h\gets$\Call{Damped\_Jacobi}{$\Amat^h,\fvec^h,\vvec^h$} 
    \ElsIf{using Gauss-Seidel}
        \State $\vvec^h\gets$\Call{Gauss-Seidel}{$\Amat^h,\fvec^h,\vvec^h$}
    \EndIf
    \State $k_\text{it}\gets k_\text{it}+1$
\EndWhile
\State \Return $\vvec^h$
\EndFunction
\Statex \Comment main function of the (damped) Jacobi and Gauss-Seidel solvers

\Function{Damped\_Jacobi}{$\Amat^h,\fvec^h,\vvec^h$}
\State $\vvec^h\gets(1-\omega)\vvec^h-\omega(\Dmat^h)\inv(\Lmat^h+\Umat^h)\vvec^h+\omega(\Dmat^h)\inv\fvec^h$
\State \Return $\vvec^h$
\EndFunction
\Statex \Comment one iteration of the damped Jacobi method

\Function{Gauss-Seidel}{$\Amat^h,\fvec^h,\vvec^h$}
\State $\vvec^h\gets-(\Lmat^h+\Dmat^h)\inv\Umat^h\vvec^h+(\Lmat^h+\Dmat^h)\inv\fvec^h$
\State \Return $\vvec^h$
\EndFunction
\Statex \Comment one iteration of the Gauss-Seidel method

\end{algorithmic}
\end{algorithm}

Based on an eigenvalue analysis of $\Gmat_\omega$ in Eq.~\ref{eqn:error_amplification}, we analyze the convergence properties of the damped Jacobi method for different modes in the error $\evec\curstep$ (i.e., eigenvectors of $\Gmat_\omega$).
We expand the error $\evec\curstep$ as
\begin{align}
    \label{eqn:error_decompose}
    \evec\curstep=\sum_{j=1}^n e_j\curstep\phivec_j= \sum_{j=1}^n \lambda_j^n e_j\initstep\phivec_j,
\end{align}
where $\{\lambda_j\}_{j=1}^n$, $\{\phivec_j\}_{j=1}^n$ are the eigenvalues and eigenvectors of $\Gmat_\omega$, respectively; $e_j\curstep$ is the coefficient of mode $\phivec_j$ in the error $\evec\curstep$. For convenience, we sort $\{\lambda_j\}_{j=1}^n$ and $\{\phivec_j\}_{j=1}^n$ by the frequencies of $\phivec_j$ in ascending order, so that $\phivec_1$ corresponds to the lowest-frequency mode and $\phivec_n$ corresponds to the highest-frequency mode.
According to Eq.~\ref{eqn:error_decompose}, the method converges if and only if $\max_{j}|\lambda_j|<1$, and the convergence rate is dominated by $\max_{j}|\lambda_j|$. For the Poisson equation defined in Eq.~\ref{eqn:poisson}
, the damped Jacobi solver yields $\max_{j}|\lambda_j|=|\lambda_1|<1$, making the iterative solution converge. For the Helmholtz equation defined in Eq.~\ref{eqn:helmholtz} with sufficiently large $k(x)$, $\max_{j}|\lambda_j|=|\lambda_1|>1$, making the iterative solution diverge.

\subsubsection{Multigrid Methods}
\label{SI_sub:MG}

To address the issue of slow convergence, multigrid methods have been proposed and widely applied in recent decades. The essence is to update the approximate solution through nested iterations in multiple grids. For a large linear system, one maps the problem into a coarser grid in a recursive approach, through which the low-frequency modes in the original grid becomes relatively high-frequency for the coarse grid, hence converging significantly faster than solving the problem using the original grid. A typical multigrid V-cycle algorithm for solving $\Amat^h\uvec^h=\fvec^h$ with (damped) Jacobi relaxation is summarized in Algorithm~\ref{alg:MG}.

\begin{algorithm}[h!]
\caption{Multigrid V-Cycle Algorithm with (Damped) Jacobi Relaxation}
\label{alg:MG}
\begin{algorithmic}

\Function{Multigrid\_Solver}{$\Amat^h,\fvec^h$}
\State $\vvec^h\gets\zerovec^h$ \Comment{initial guess of the solution}
\State $k_\text{v}\gets 1$
\While{$k_\text{v}\leq n_\text{v}$ and not converged}
    \State $\rvec^h\gets\fvec^h-\Amat^h\vvec^h$
    \State $\vvec^h\gets$\Call{V\_cycle}{$\Amat^h,\fvec^h,\vvec^h$} 
    \State $k_\text{v}\gets k_\text{v}+1$
\EndWhile
\State \Return $\vvec^h$
\EndFunction
\Statex \Comment main function of the multigrid V-cycle algorithm

\Function{V\_cycle}{$\Amat^h,\fvec^h,\vvec^h$}
\State Loop $n_\text{it}$ times $\vvec^h=$\Call{Damped\_Jacobi}{$\Amat^h,\fvec^h,\vvec^h$} \Comment see Alg.~\ref{alg:JGS}
\If{$\Omega^h$ is not the coarsest grid}
\State $\Amat^{2h},\fvec^{2h}\gets$\Call{restriction}{$\Amat^h$}, \Call{restriction}{$\fvec^h-\Amat^h\vvec^h$}
    \State $\vvec^{2h}\gets$\Call{v\_cycle}{$\Amat^{2h},\fvec^{2h},\zerovec^{2h}$}
    \State $\vvec^h\gets\vvec^h+$\Call{prolongation}{$\vvec^{2h}$}
\EndIf
\State Loop $n_\text{it}$ times $\vvec^h=$\Call{Damped\_Jacobi}{$\Amat^h,\fvec^h,\vvec^h$}
\State \Return $\vvec^h$
\EndFunction
\Statex \Comment one V-Cycle
\end{algorithmic}
\end{algorithm}

Multigrid methods have been shown to work for problems such as the Poisson equation, where the left-hand-side matrix $\Amat^h$ is positive definite. However, there exists a wide class of problems where $\Amat$ is not positive definite. As a result, yielding $\max_j|\lambda_j|>1$ so that the multigrid method based on Jacobi relaxations diverges.

\subsubsection{General Summary of DeepONet}
\label{SI:deeponet_general}

\begin{figure}[H]
  \centering
  \includegraphics[width=0.5\textwidth]{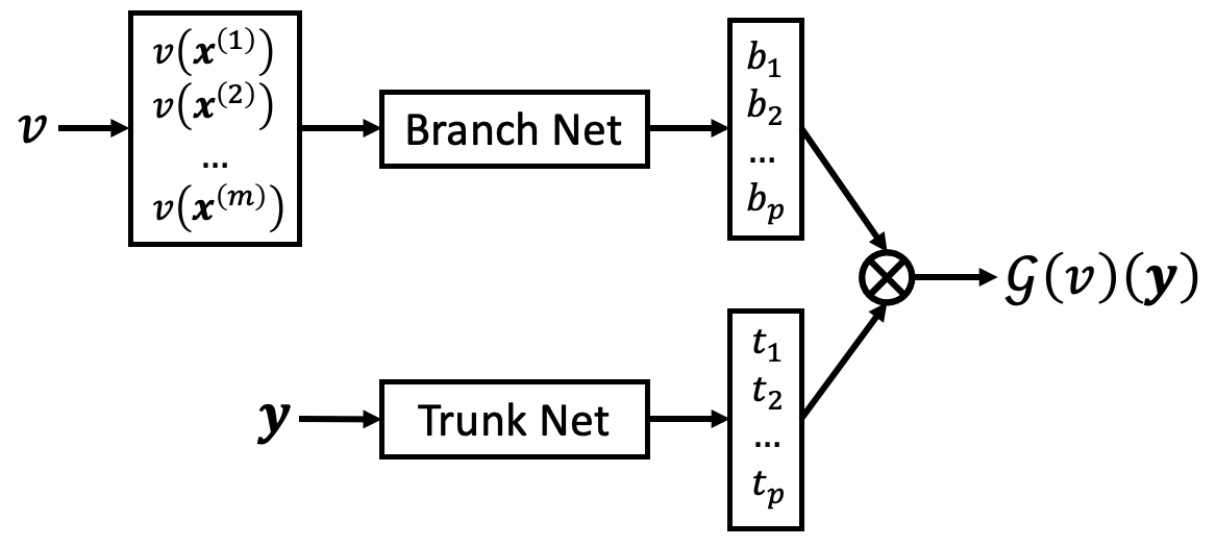}
\caption{\textbf{Architecture of a DeepONet.} The DeepONet learns the operator $\Gopt$ that maps the input function $v$ into the output function $u:=\Gopt(v)$. The branch network takes as input $[v(\xvec^{(1)}),v(\xvec^{(2)}),...,v(\xvec^{(m)})]\tran$, i.e., the evaluation of function $v$ at $\xvec^{(1)},\xvec^{(2)},...,\xvec^{(m)}$; the trunk network takes $\yvec$ as input. Branch network and trunk network yield outputs $[b^{(1)},b^{(2)},...,b^{(p)}]\tran$ and $[t^{(1)},t^{(2)},...,t^{(p)}]\tran$, respectively. The inner product of these two vectors (plus a trainable bias) is the approximation of $\Gopt(v)(\yvec)$ as the final output of the DeepONet.}
  \label{fig:deeponet}
\end{figure}

DeepONet approximates operators, i.e., the mapping from functions to functions. We define the target nonlinear operator as $\Gopt: v \mapsto u$, where $v: \xvec \mapsto v(\xvec) \in \Rdomain$ and $u: \yvec \mapsto u(\yvec) \in \Rdomain$ are the input and output functions that are defined for $\xvec \in D \subset \Rdomain^d$ and $\yvec \in D' \subset \Rdomain^{d'}$, respectively. The architecture of a DeepONet is illustrated in Fig.~\ref{fig:deeponet}. The DeepONet comprises of two subnetworks, namely, a branch network and a trunk network. To approximate $u(\yvec)=\Gopt(v)(\yvec)$ for $\yvec \in D'$, the branch network takes as input $[v(\xvec^{(1)}),v(\xvec^{(2)}),...,v(\xvec^{(m)})]\tran$, i.e., the evaluation of function $v$ at $\xvec^{(1)},\xvec^{(2)},...,\xvec^{(m)}\in D$; the trunk network takes $\yvec \in D'$ as input. Branch network and trunk network yield outputs $[b^{(1)},b^{(2)},...,b^{(p)}]\tran$ and $[t^{(1)},t^{(2)},...,t^{(p)}]\tran$, respectively. The inner product of these two vectors is the approximation of $\Gopt(v)(\yvec)$ as the final output of the DeepONet. Complete details regarding DeepONet are explained in the original publication~\cite{lu2021deeponet}.

\subsubsection{Pseudo-code of HINTS algorithms}
\label{SI:technical}

\begin{algorithm}[H]
\caption{HINTS-Jacobi}
\label{alg:HINTS-Jacobi}
\begin{algorithmic}
\Function{HINTS\_Jacobi}{$k(\xvec),f(\xvec),\text{PDE}$}
\State $\Amat^h,\fvec^h$ = \Call{discretize}{$k(\xvec),f(\xvec),\text{PDE}$} \Comment{discretization; to solve $\vvec^h=(\Amat^h)^{-1}\fvec^h$}
\State $\vvec^h\gets\zerovec^h$ \Comment{initial guess of the solution}
\State $k_\text{it}\gets 1$
\While{$k_\text{it}\leq n_\text{it}$ and not converged}
        \State $\rvec^h\gets\fvec^h-\Amat^h\vvec^h$
    \If{$k_\text{it}\text{ mod }n_{\text{r}} = 0$} \Comment{condition for invoking DeepONet}
        \State $r(\xvec)=$\Call{reverse\_discretize}{$\rvec^h$}  \Comment{see Section S3E}
        \State $\delta \vvec^h$ ($=\delta v(\xvec)$)$\gets$\Call{DeepONet}{$k(\xvec),r(\xvec)$}
        \State $\vvec^h\gets \vvec+\delta\vvec^h$
    \Else
        \State $\vvec^h\gets$\Call{Damped\_Jacobi}{$\Amat^h,\fvec^h,\vvec^h$} \Comment see Alg.~\ref{alg:JGS}
        \State (or equivalently: $\vvec^h\gets\vvec^h+$\Call{Damped\_Jacobi}{$\Amat^h,\rvec^h,\zerovec^h$})
    \EndIf
    \State $k_\text{it}\gets k_\text{it}+1$
\EndWhile
\State \Return $\vvec^h$ ($=v(\xvec)$)
\EndFunction
\end{algorithmic}
\end{algorithm}

\begin{algorithm}[H]
\caption{HINTS-MG with V-cycle and Jacobi Relaxation}
\label{alg:HINTS-MG}
\begin{algorithmic}
\Function{HINTS\_MG}{$k(\xvec),f(\xvec),\text{PDE}$}
\State $\Amat^h,\fvec^h$ = \Call{discretize}{$k(\xvec),f(\xvec),\text{PDE}$}
\State $\vvec^h\gets\zerovec^h$
\State $k_\text{v}\gets 1$
\While{$k_\text{v}\leq n_\text{v}$ and not converged}
    \State $\vvec^h=$\Call{HINTS\_V\_cycle}{$\Amat^h,\fvec^h,\vvec^h,k(\xvec)$}
\EndWhile
\EndFunction
\Statex \Comment main function
\Function{HINTS\_V\_cycle}{$\Amat^h,\fvec^h,\vvec^h,k(\xvec)$}
    \State $\vvec^h=$\Call{HINTS\_Jacobi\_discrete}{$\Amat^h,\fvec^h,\vvec^h,k(\xvec)$}
    \Statex \Comment Modified Alg.~\ref{alg:HINTS-Jacobi} without the first discretization step
    \If{$\Omega^h$ is not the coarsest grid}
        \State $\Amat^{2h},\fvec^{2h}\gets$\Call{restriction}{$\Amat^h$}, \Call{restriction}{$\fvec^h-\Amat^h\vvec^h$}
        \State $\vvec^{2h}\gets$\Call{HINTS\_V\_cycle}{$\Amat^{2h},\fvec^{2h},\zerovec^{2h},k(\xvec)$}
        \State $\vvec^h\gets\vvec^h+$\Call{prolongation}{$\vvec^{2h}$}
    \EndIf
    \State $\vvec^h=$\Call{HINTS\_Jacobi\_discrete}{$\Amat^h,\fvec^h,\vvec^h,k(\xvec)$}
    \State $k_\text{v}\gets k_\text{v}+1$
\State \Return $\vvec^h$
\EndFunction
\Statex \Comment a single V-Cycle in HINTS-MG
\end{algorithmic}
\end{algorithm}

\subsubsection{Discretization of Differential Equation}

\begin{figure}[H]
  \centering
  \includegraphics[width=0.95\textwidth]{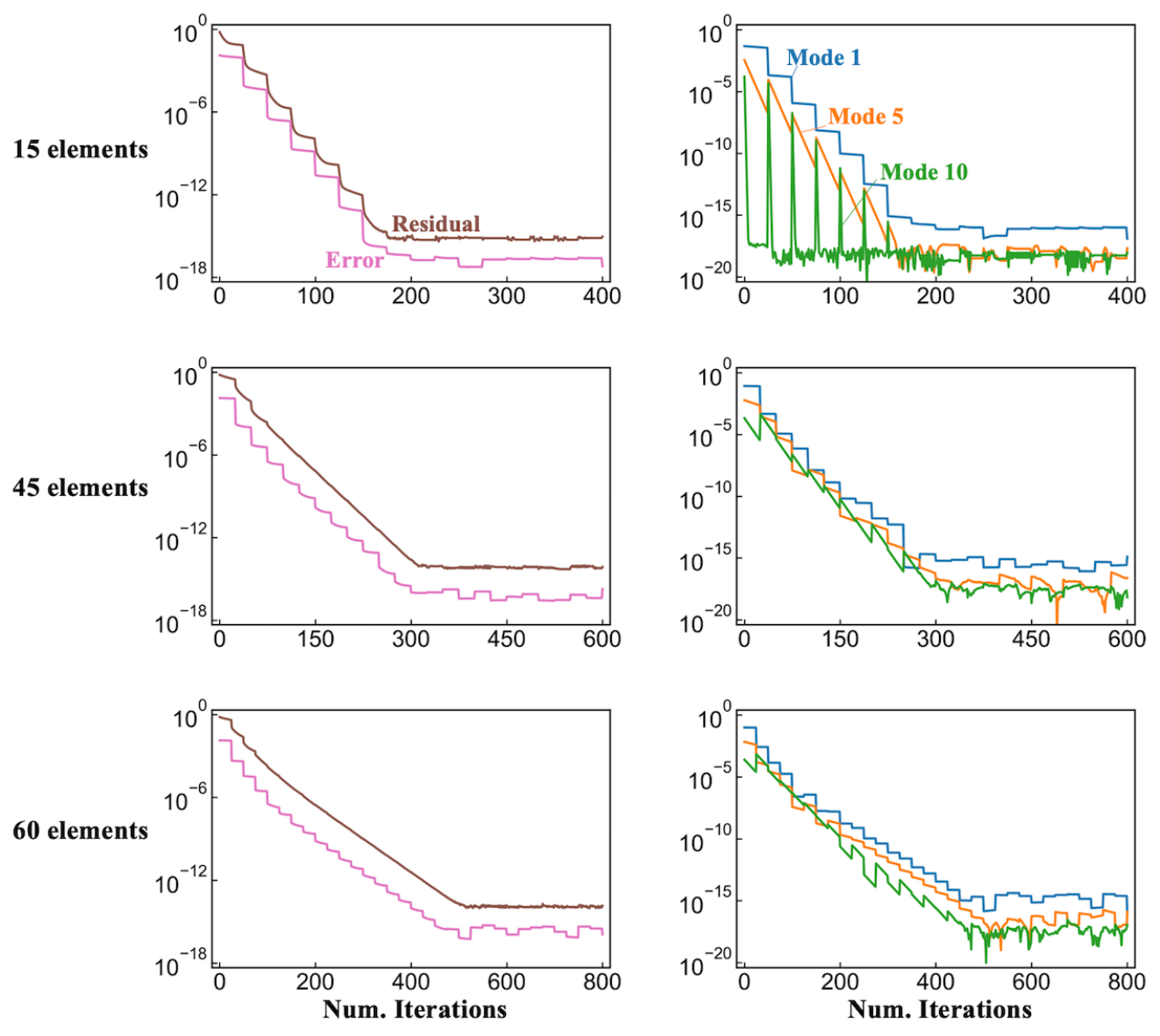}
\caption{\textbf{Results of 1D Poisson Equation with Mismatching Discretization.} We use the identical DeepONet to the one in the main text (trained with $n_\text{D}=30$) to solve equations discretized with different numbers of intervals $n=15,45,60$. The first column and the second column shows the histories of the norms of residual and error (first column) and the mode-wise error (second column), respectively.}
  \label{fig:results_1P_NoE}
\end{figure}

We use the finite element method to discretize the 1D Poisson equation (Eq.~\ref{eqn:poisson}). The computational domain $\Omega=(0,1)$ is partitioned into $n=30$ linear elements, except in Fig.~\ref{fig:results_1P_NoE} where $n=15,45,60$.

For the 2D Poisson equation, we use the finite element method to discretize the domain into linear triangular elements (see the first panel in Fig.~\ref{fig:results_HD}(A)). The mesh is generated using MATLAB. The side lengths for elements far from the corner $(0,0)$ are no larger than $0.05$, while those for elements near the corner are $0.005$. In total, there are 450 nodes and 802 elements in the domain.

For the $d$-dimensional ($d\in\{1,2,3\})$) Helmholtz equation (Eq.~\ref{eqn:helmholtz}), we discretize the domain $\Omega=(0,1)^d$ into $n=30,32^2,16^3$ Cartesian grids, respectively.
\revv{The Helmholtz-annular cylinder example is constructed by discretizing (Eq.~\ref{eqn:helmholtz}) using the finite element method and unstructured mesh encapsulating the annular cylinder, depicted on Fig.~\ref{fig:helm3D_cylinder_example}(A). We impose zero Dirichlet BC only on the top of the cylinder.}

\subsubsection{DeepONet Architecture}
\label{SI:don_architecture}

DeepONet takes the input function $k(\xvec)$ and $f(\xvec)$ in the forms of $[k(\xvec_0),k(\xvec_1),...,k(\xvec_{n_{\text{D}}})]\tran$ and $[f(\xvec_0),f(\xvec_1),...,f(\xvec_{n_{\text{D}}})]\tran$, respectively, where $\{\xvec_i\}_{i=0}^{n_{\text{D}}}$ are uniform grid points in $\Omega$. For 1D cases, the size of the branch network is $[62, 60, 60, 60]$, where we concatenate $k(x)$ and $f(x)$, resulting in the input dimension $2(n_{\text{D}}+1)=62$; the dimension of the trunk network is $[1, 60, 60, 60]$. We train the DeepONet for $10,000$ epochs with an initial learning rate $10^{-3}$, a decay rate of $50\%$ per $5,000$ epochs, and a mini-batch size of $500$. See the training results for the 1D Poisson equation in Fig.~\ref{fig:deeponet_training}. 
\textcolor{red}{ Creating the data to train the 1D DeepONets in the paper takes a few seconds, and training for 10,000 epochs takes about a minute. A total offline training time for a 1D HINTS is roughly 2 minutes.}

For the 2D cases, the branch net is a combination of a CNN (input dimension $31\times 31$, number of channels $[2, 40, 60, 100, 180]$, kernel size $3\times 3$, stride 2, the number of channels of the input $2$ comes from the concatenation of $k(x)$ and $f(x)$) and dense NN (dimension $[180,80,80]$). The dimension of the trunk network is [2, 80, 80, 80]. We train the DeepONet for $25,000$ epochs with a fixed learning rate $10^{-4}$, and a mini-batch size of 10,000.

For the 3D Helmholtz \oldtext{case}\revv{examples} \oldtext{in Fig.~\ref{fig:results_HD}}, the branch net is a combination of a 3D CNN (input dimension $31\times 31\times 31$, number of channels $[2,40,60,100]$, kernel size $3\times 3\times 3$, stride 2) and dense NN (dimension $[100,80,80]$). The dimension of the trunk network is $[3, 80, 80, 80]$. $k(x,y)$ and $f(x,y)$, each with size $14\times 14\times 14$, are concatenated as two input channels in the branch network. We train the DeepONet for $25,000$ epochs with a learning rate $10^{-4}$, and a mini-batch size of 10,000. 
\textcolor{red}{ Creating the training data for the 3D DeepONet took $82$ minutes on a CPU node, and training took $160$ minutes on a NVIDIA A6000 GPU. This network serves the largest scale example we present in the manuscript, and gives an estimate for the training duration. We stress that the DeepONet smooths the high frequency modes, and thus does not require fine discretization. Moreover, training can be accomplished with fewer epochs; more epochs are usually required to capture the smaller features, corresponding to the higher frequencies in the data. While parallel training, a smaller network, better hardware, etc., could reduce training time, a $\sim 4$ hours total offline training is manageable.}

\textcolor{red}{For the large-scale 3D Helmholtz example\revv{s} in Fig.~\ref{fig:results_3H_GS}(B) and \revv{Fig.~\ref{fig:helm3D_cylinder_example}(A)}, we use the same DeepONet mentioned in the last paragraph. We switch most of the operations in the code for this case to work on sparse matrices, to be able to hold large matrices (up to $10^6\times 10^6$ and over) in the memory. We ran the code using hardware equipped with the Nvidia A6000 GPU. 
We note that the DeepONet was not trained for example depicted in Fig.~\ref{fig:results_3H_GS}(B)  with the small box cutout, but in inference we did use it, showing the generalization capabilities of the method.}

\textcolor{red}{For the large-scale 3D Poisson example in Fig.~\ref{fig:results_1P_gen_corlen}(B), we trained the DeepONet we use $k(x, y)$ and $f(x, y)$ of size $32\times32\times32$, concatenated the same as the input channels. The DeepONet architecture is the same with branch dense dimension of $[100, 256, 256, 256]$, and trunk dimension of $[3, 180, 180, 180]$. We train the DeepONet with a mini-batch size of $8,000$. The rest is similar to the foregoing large-scale 3D Helmholtz case.}

For all the cases, the activation function is the Rectified Linear Unit (ReLU) for branch networks (excluding the last layer that has no activation function) and Hyperbolic Tangent (tanh) for trunk networks (including the last layer). The loss function is defined by
\begin{align}
    \label{eqn:loss}
    \mathcal{L}=\frac{1}{N(n_{\text{D}}+1)}\sum_{j=1}^N\sum_{i=0}^{n_\text{D}}\frac{(\hat{u}_{(j)}(\xvec_i)-u_{(j)}(\xvec_i))^2}{\varepsilon+|u_{(j)}(x)|^\alpha},
\end{align}
where $N$ is the size of training data, $(\cdot)_{(j)}$ refers to quantities from the $j$th training sample, $\hat{(\cdot)}$ refers to the prediction, $\alpha$ is a parameter, and $\varepsilon$ is a small number for preventing zero denominator when $\alpha\neq 0$. For the 1D Poisson equation, $\alpha=0$, making the loss function the mean squared error. For the 1D Helmholtz equation, $\alpha=1$. For other cases, $\alpha=2$. When $u(x)$ changes dramatically across different samples, a positive value of $\alpha$ yields weights that are biased appropriately towards $u(x)$ with smaller magnitudes.

\subsubsection{Conversion of Residuals}
Standard numerical solvers deal with discretized, matrix- and vector-form representation of the system, such as $\Amat^h$, $\fvec^h$ or $\rvec^h$, and $\uvec^h$. On the other hand, the DeepONet solver deals with the function-form representation of the system, such as $k(\xvec)$, $f(\xvec)$ or $r(\xvec)$, and $u(\xvec)$ (more specifically, the discrete evaluation of these functions).

We now explain in detail how to bridge the gap between these two representations. First, $u(\xvec)$ and $\uvec^h$ are essentially the same quantities (up to $n\neq n_\text{D}$ which simply needs an interpolation). The only issue is in terms of the conversion between $\fvec^h$ ($\rvec^h$) and $f(\xvec)$ ($r(\xvec)$). For all cases except the 2D Poisson equation in an L-shaped domain, the grids are uniform, so that $\fvec^h$ and $f(\xvec)$ are possibly different only in terms of a constant of the grid size. We always absorb this constant into $\Amat^h$, so that $\fvec^h$ and $f(\xvec)$ are equivalent (similar to $u(\xvec)$ and $\uvec^h$). For the case of the 2D Poisson equation in an L-shaped domain, due to the nonuniform meshes, we need to explicitly calculate $f(\xvec)$ from $\fvec^h$. We approximately calculate $f(\xvec_i)$, where $\xvec_i$ is node $i$, by
\begin{align}
    f(\xvec_i)= \frac{3\text{f}^h_i}{\sum_j A_j},
\end{align}
where $\text{f}^h_i$ is the component $i$ of $\fvec^h$, $A_j$ is the area of element $j$, and the summation is calculated for all elements that include node $i$. While this calculation is approximate, it is only used for the DeepONet inference that essentially provides an approximate solution. Hence, such approximation does not influence the overall performance of HINTS.

\subsection{Supplementary Tables}

\subsubsection{Data Generation and DeepONet Training}
\label{SI:datageneration}

The input data $(k(\xvec),f(\xvec))$ for training and testing cases are generated using GRFs with means $\mathbb{E}k(\xvec)=k_0$ ($>0$), $\mathbb{E}f(\xvec)=0$ and covariances
\begin{align}
    \label{eqn:cov_k}
    \Cov\big(k(\xvec_1),k(\xvec_2)\big)&=\sigma_k^2\exp\Big(-\frac{||\xvec_1-\xvec_2||^2}{2l_k^2}\Big);\\
    \label{eqn:cov_f}
    \Cov\big(f(\xvec_1),f(\xvec_2)\big)&=\sigma_f^2\exp\Big(-\frac{||\xvec_1-\xvec_2||^2}{2l_f^2}\Big).
\end{align}
In addition, we only consider $k(\xvec)$ satisfying $k(\xvec)>k_\text{min}$. The choices of $k_0$, $k_\text{min}$, $\sigma_k$, $l_k$, $\sigma_f$ $l_f$ are listed in Table~\ref{table:GRF}.

\begin{figure}[H]
  \centering
\includegraphics[width=\textwidth]{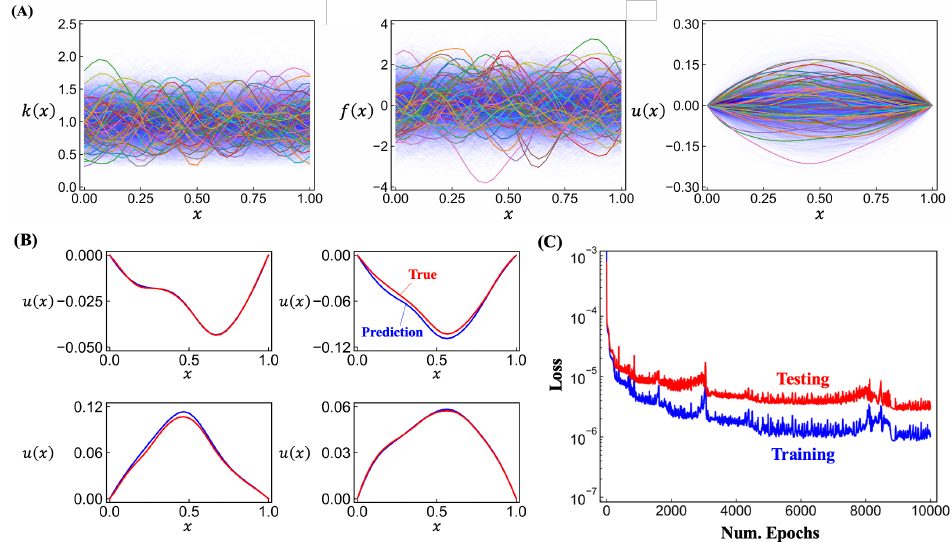}    
\caption{
\textbf{DeepONet Training Data for 1D Poisson Equation.}
(A) An example of functions $k(x), f(x)$, generated by the GRFs and associated solutions $u(x)$.
(B) Comparison of the predicted and true solutions $u(x)$ for four typical test cases with different $k(x)$ and $f(x)$. 
(C) Loss function for training and testing.}
  \label{fig:deeponet_training}
\end{figure}

\begin{table}[H]
\small
\centering
\begin{tabular}{c c c c c c c} 
 \hline
 \hline
  Case & $k_0$ & $k_\text{min}$ & $\sigma_k$ & $l_k$ & $\sigma_f$ & $l_f$ \\ 
 \hline
    1D Poisson & $1.0$ & $0.3$ & $0.3$ & $0.1$ & $1.0$ & $0.1$  \\
    1D Helmholtz & $8.0$ & $3.0$ & $2.0$ & $0.2$ & $1.0$ & $0.1$ \\
    2D Poisson (L-shaped) & $1.0$ & $0.3$ & $0.3$ & $0.1$ & $1.0$ & $0.1$ \\
    2D Helmholtz & $6.0$ & $3.0$ & $0.5$ & $0.3$ & $1.0$ & $0.1$ \\
    3D Helmholtz* & $6.0$ & $3.0$ & $0.5$ & $0.5$ & $1.0$ & $0.1$ \\
    \revv{3D Helmholtz (Annular cylinder)} & \revv{$1.0, 3.0, 6.0$} & \revv{$0.0$} & \revv{$0.5$} & \revv{$0.5$} & \revv{$1.0$} & \revv{$0.1$}\\
 \hline
 \hline
\end{tabular}
\caption{\textbf{Parameters for Generating Training Data.} *The parameters shown are for test data. For the training data, $k_0=6.0$, $k_\text{min}=0.3$, $\sigma_k=0.2$, $l_k=0.1$, $\sigma_f=0.1$, $l_f=0.2$. The dataset for Annular cylinder example has been uploaded to Zenodo repository and can be freely accessed at \href{https://doi.org/10.5281/zenodo.10904349}{https://doi.org/10.5281/zenodo.10904349}}
\label{table:GRF}
\end{table}
After $(k(\xvec)$ and $f(\xvec))$ are generated, we solve $u(\xvec)$ using standard numerical solvers. 
As an example, the generated training data for the 1D Poisson equation is shown in Fig.~\ref{fig:deeponet_training}(A). 
For all the Helmholtz cases, these choices of $k(x)$ ensure the indefiniteness of the system.

\subsubsection{Summary of notation}
\label{SI:notations}
Table~\ref{table:notations_a} summarizes the notations used in the manuscript.

{
\begin{longtable}{|p{2cm}|p{10.8cm}|}
\hline
    $x$ ($\xvec$) & Spatial coordinates for one- (multi-) dimensional problems \\
    $\Amat^h$ & Matrix in the linear system $\in\mathbb{R}^{n\times n}$; with discretization $h$) \\
    $\Lmat^h$ & Strictly lower triangular part of $\Amat^h$ \\
    $\Dmat^h$ & Diagonal part of $\Amat^h$ \\
    $\Umat^h$ & Strictly upper triangular part of $\Amat^h$ \\
    $\Mmat^h, \Nmat^h$ & Two matrices split from $\Amat^h$ (see Eqs.~\ref{eqn:iterative_general}-~\ref{eqn:iterative_MN_GS}) \\
    $\fvec^h$ & Right-hand-side vector in the linear system  ($\in\mathbb{R}^{n\times 1}$; with discretization $h$) \\
    $\rvec^h$ & Residual vector \\
    $\evec^h$ & Error vector \\
    $\uvec^h$ & True solution of linear system $\Amat^h\uvec^h=\fvec^h$\\
    $\vvec^h$ & Approximate solution of $\uvec^h$\\
    $\delta\vvec^h$ & Correction of $\vvec^h$ in an iteration\\
    $\phivec_i^h$ & Eigenmode $i$ of $\uvec^h$ \\
    $(\cdot)^{h}$ & Matrices/Vectors corresponding to the grid $\Omega^{h}$ \\
    $(\cdot)^{2h}$ & Matrices/Vectors corresponding to the coarser grid $\Omega^{2h}$ (with respect to $\Omega^{h}$) \\
    $k=k(\xvec)$ & function parameterizing the differential equation (see Eq.~\ref{eqn:diffeqn_general_inside})\\
    $f=f(\xvec)$ & function parameterizing the differentiql equation (see Eq.~\ref{eqn:diffeqn_general_inside}) \\
    $r=r(\xvec)$ & residual function\\
    $u=u(\xvec)$ & True solution (function) of the differential equation (see Eqs.~\ref{eqn:diffeqn_general_inside}-\ref{eqn:diffeqn_general_bc})\\
    $v=v(\xvec)$ & Approximate solution (function) of $u(\xvec)$\\
    $\phi_i=\phi_i(\xvec)$ & Eigenmode $i$ of $u(\xvec)$ \\ 
    $l_k$, $l_f$ & Correlation lengths of the GRF for generating $k(x)$ and $f(x)$ \\
    $N$ & Number of training cases \\
    $n$ & Number of grids/elements/sub-intervals in the numerical problem \\
    $n_\text{D}$ & Number of grids/elements/sub-intervals for DeepONet training \\
    $k_{\text{it}}$ & The index of the current iteration in the numerical solver\\
    $k_{\text{v}}$ & The index of the current V-cycle in the multigrid solver\\
    $(\cdot)^{(k_{\text{it}})}$ & Matrices/Vectors in $k_{\text{it}}$-th iteration \\
    $n_{\text{it}}$ & Number of total iterations in the numerical solver\\
    $n_{\text{r}}$ & Reciprocal of the proportion of DeepONet in HINTS (i.e., DeepONet-to-numerical ratio $1:(n_{\text{r}}-1)$\\
    $n_{\text{v}}$ & Number of V-cycles in multigrid methods\\
    $n_{\text{rl}}$ & Number of relaxation steps in each grid level of multigrid methods \\
    $n_{\text{cut}}$ & Number of modes that the DeepONet is able to predict accurately \\
    $\mu$ & Convergence rate (see Eq.~\ref{eqn:conv_rate}) \\
    $\Gopt$ & Operator that the DeepONet seeks to learn (see Eq.~\ref{eqn:deeponet})\\
    $\Lopt_\xvec$ & Linear differential operator in Eq.~\ref{eqn:diffeqn_general_inside}\\
    $\Bopt_\xvec$ & Boundary operator in Eq.~\ref{eqn:diffeqn_general_bc}\\
    $\Omega$ & computational domain of the differential equation\\
    $\Omega^{\hDON}$ & Discretized $\Omega$ for DeepONet \\
    $\Omega^h$ & Discretized $\Omega$ for the numerical problem\\
    $\omega$ & Damping parameter\\
    $\alpha$ & Parameter in the loss  function (see Eq.~\ref{eqn:loss})\\
    $\hat{(\cdot)}$ & Predicted value from DeepONet\\
    $(\cdot)_{(j)}$ & Quantity from the $j$th training sample\\ \hline
\caption{Summary of mathematical notations utilized in manuscript.}
\label{table:notations_a}
\end{longtable}
}

\subsection{Supplementary Discussion}
\subsubsection{Eigenmode Analysis of HINTS}
\label{SI:eigenAnalysis}
\begin{figure}[h!]
  \centering
\includegraphics[width=\textwidth]{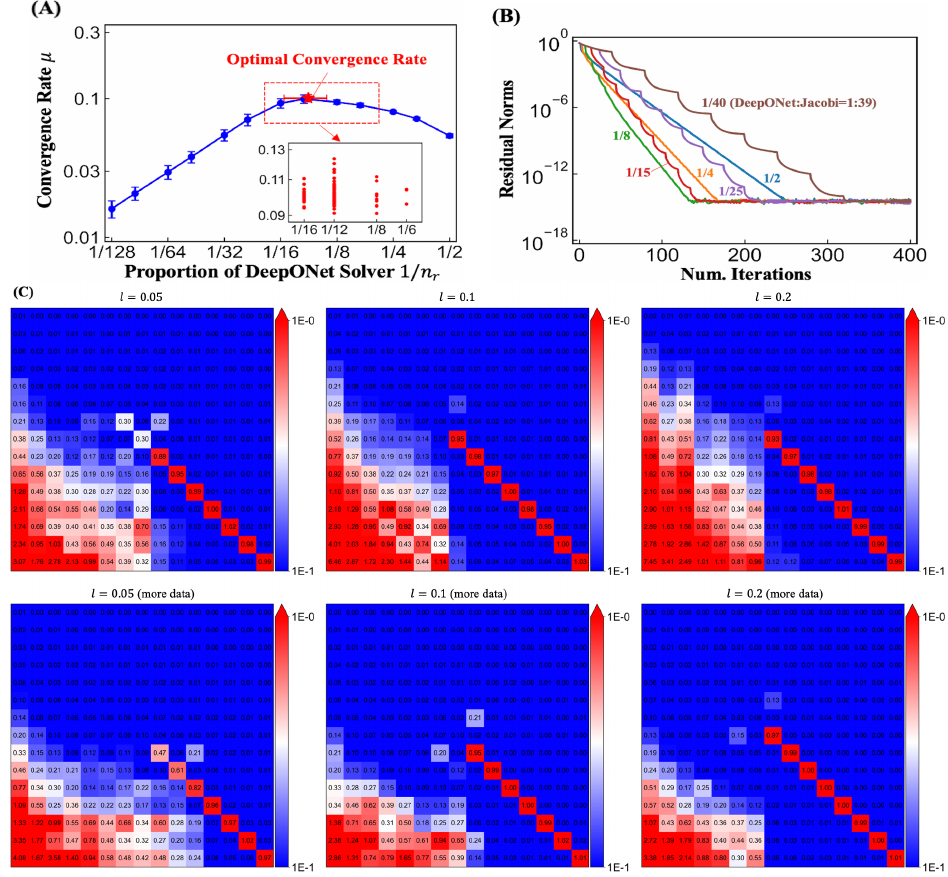}    
\caption{\textbf{Impact of different DeepONet proportions and residual sampling strategies on the convergence of HINTS.}
The results are reported for the test case $(k,f)$ used in Fig.~\ref{fig:results_1P}(A). 
(A) Convergence rates $\mu$ for different proportions of DeepONet ($1/n_\text{r}$) in HINTS-Jacobi. 
The mean value and standard deviation of $\mu$ for $100$ test cases is shown in blue. 
The optimal convergence rate is marked by red. 
(B) Convergence histories of residual norms for different proportions of DeepONet in HINTS-Jacobi.
(C) Mode-wise performance of the trained DeepONet.
Row $j$ ($j\in\{1,2,...,15\}$) represents the test of DeepONet with input $f(x)$ corresponding to true solution $u(x)=\phi_j(x)$, i.e., only the eigenmode $j$. 
The mode-wise errors of the DeepONet prediction/output are shown in different columns of matrix.
The first/second row corresponds to 5,000/10,000 training cases. 
The columns depict different results for the training data generated from GRFs with correlation lengths $l=0.05, 0.1, 0.2$, respectively.}
  \label{fig:results_1P_analysis_eig}
\end{figure}

\begin{figure}[H]
  \centering
 \includegraphics[width=\textwidth]{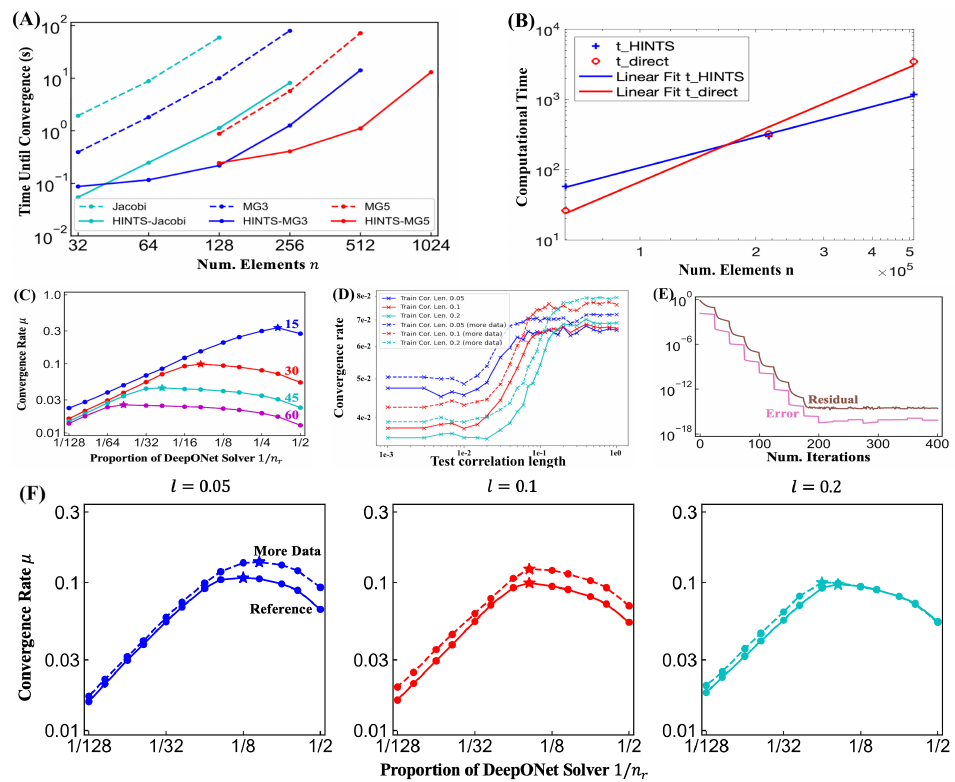}  
%
%
%
%
\caption{\textbf{Extended results for the 1D Poisson equation.}
(A) Computational cost of multigrid (MG) and HINTS-MG methods.
All MG methods use damped Jacobi relaxation. MG3/MG5 and HINTS-MG3/MG5 use three/five grid levels.
(B) Computational time as a function of system size (3D problem) for direct solver and HINTS-Jacobi.
(C) The convergence rate of HINTS-Jacobi as a function of different DeepONet proportions. 
The four lines correspond to numbers of sub-intervals $n = 15, 30, 45, 60$. 
(D) The convergence rate of HINTS-Jacobi with respect to out-of-distribution input functions $k(x)$ and $f(x)$ with different correlation lengths. Different curves refer to different training setups, including different correlation lengths and volumes of the training data.
(E) Convergence of HINTS-GS.
(F) Convergence rate of HINTS-Jacobi with respect to different proportion of the DeepONet and different residual sampling strategies. 
We consider different correlation lengths ($l=0.05,0.1,0.2$) in the GRF for generating the training data, as well as different size of training data (5,000 or 10,000 cases).
}
\label{fig:results_1P_gen_corlen}
\end{figure}

From the mode-wise errors in Fig.~\ref{fig:results_1P}(C), we have shown that DeepONet performs well for the low-frequency modes in the solution. This characteristic is associated with the fact that the DeepONet training data are generated from a GRF with a certain correlation length, which mainly contains frequencies within a certain range. Previous studies on DeepONet have found that DeepONet generalizes well towards frequencies lower than the training data but not for higher frequencies~\cite{don_bubble,don_coupling}, which is closely related to spectral bias of neural networks~\cite{rahaman2019spectral} that we will further address in Discussion. Hence, we expect that the trained DeepONet can predict well solutions of mode 1, 2, ..., $n_\text{cut}$, where $n_\text{cut}$ is a cutoff frequency. 
To further demonstrate this, we examine the mode-wise performance of DeepONet in Fig.~\ref{fig:results_1P_analysis_eig}(C). 
Row $j$ ($j\in\{1,2,...,15\}$) represents the test of DeepONet with the input $f(x)$ corresponding to $u(x)=\phi_j(x)$, i.e., the solution with eigenmode $j$ only. The mode-wise errors of the DeepONet prediction/output are shown in different columns. The displayed error values are the geometric mean of 100 test cases of $k(x)$. We note that $n_\text{cut}=7$ for the current setup. For tests on the first seven modes (rows 1 to 7), DeepONet uniformly decreases errors in all modes. For tests starting from the eighth mode, DeepONet is unable to decrease the error of the input mode (shown in diagonal blocks), and even produces errors in modes 1 to 7. 
We point out that this $n_\text{cut}$ relies on the performance of DeepONet, which is influenced by diverse factors in its training process, such as size and correlation length of training data. 
We include a systematic study on these factors in Figures~\ref{fig:results_1P_gen_corlen}(C-F).

We found that a mixture of solvers helps accelerates the solution. A natural question is: what is the optimal proportion of DeepONet in HINTS? 
Intuitively, based on the preceding analysis, too large a proportion of DeepONet does not allow the Jacobi solver to cause a sufficient decay of high-frequency errors, so that the DeepONet receives inputs contaminated by high-frequency modes beyond its cutoff frequency. 
On the other hand, a too small proportion of DeepONet causes excessive Jacobi iterations that suffer from slow convergence for low-frequency modes. 
In Fig.~\ref{fig:results_1P_analysis_eig}(B), we show the residual norms for six different proportions of DeepONet, ranging from $1/40$ to $1/2$, for the test case in Fig.~\ref{fig:results_1P}(A). 
We observe non-monotonic convergence rates of convergence with respect to DeepONet proportions. To better quantify this, we define the convergence rate $\mu$ of an iterative solver by
\begin{align}
    \label{eqn:conv_rate}
    \mu:=-\frac{1}{k_\text{it,start}-k_\text{it,end}}\log_{10}\frac{||\rvec^{h{(k_\text{it,end})}}||}{||\rvec^{h{(k_\text{it,start})}}||},
\end{align}
where $k_\text{it,start}$ and $k_\text{it,end}$ are two snapshots between which the residual decays steadily. Intuitively, this convergence rate $\mu$ represents how many orders of magnitude of the residual are decayed in each iteration, which is proportional to the slopes of lines in figures such as Fig.~\ref{fig:results_1P_analysis_eig}(B). We show the statistics of $\mu$ for different DeepONet proportions $1/n_\text{r}$ for all 100 test cases in Fig.~\ref{fig:results_1P_analysis_eig}(A). In this particular setup, the optimal proportions range from $1/16$ to $1/6$, with the mean value being approximately $1/12$.

\subsubsection{Additional Results}
\label{SI:directresults}
This section presents additional results for the 1D Poisson equation, including those for (1) other $k(x)$ and $f(x)$, (2) HINTS-GS and (3) different discretizations for DeepONet and the numerical problem (i.e., $\Omega^{\hDON}\neq\Omega^h$, or $n_\text{D}\neq n$). 
The performance of HINTS is consistent with the results in Fig.~\ref{fig:results_1P} for both cases, which demonstrates (1) effectiveness of HINTS in its integration with the Gauss-Seidel method, and (2) the transferability/generalizability of HINTS across different discretizations of the computational domain.

\begin{figure}[H]
  \centering
  \includegraphics[width=\textwidth]{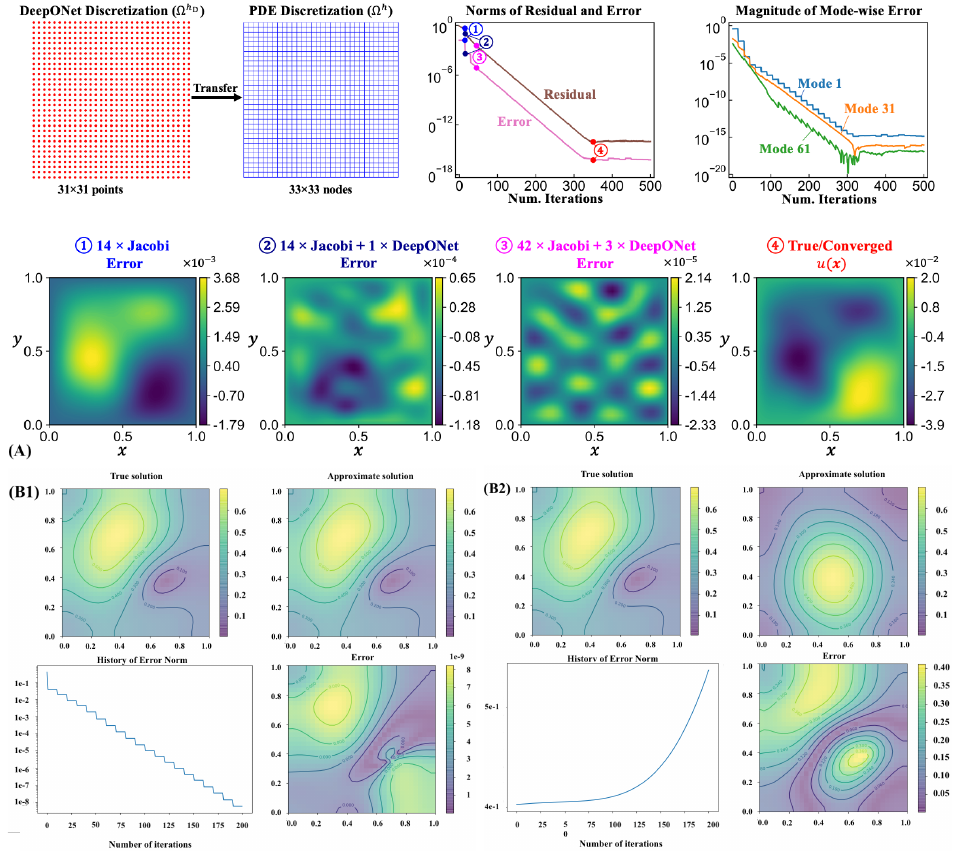}  
\caption{\textbf{Results of 2D Helmholtz Equation with HINTS-Jacobi.}
(A) Helmholtz Equation with Dirichlet boundary conditions. 
(B) Helmholtz Equation with   Sommerfeld radiation absorbing boundary condition. 
(B1) The results obtained using HINTS-Jacobi.
(B2) The results obtained using Jacobi. }
  \label{fig:results_2H_sommerfeld}
\end{figure}

\begin{figure}[H]
  \centering
\includegraphics[width=\textwidth]{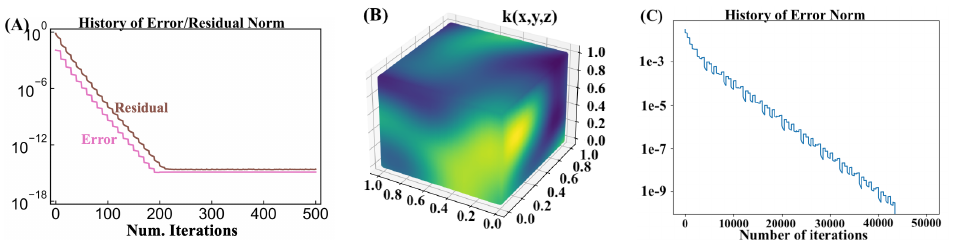}  
\caption{\textbf{Results of HINTS for 3D Helmholtz Problem.}
(A) Convergence history of HINTS-GS for Helmholtz problem with Dirichlet boundary conditions. The system is of size $16\times16\times16$, and the DeepONet is trained on examples of size $12\times12\times12$.
(B,C) A system with $\mathcal{O}(10^6)$ elements is solved using the HINTS-Jacobi. In addition, a small scatterer is added inside of the computational domain. The DeepONet for the HINTS is of size $12\times12\times12$. 
(B) An example of an input function $k$. Note that the scatterer is hidden inside and not visible in this panel. 
(C) The convergence history of the error.
}
  \label{fig:results_3H_GS}
\end{figure}

\subsubsection{Generalizability and Scalability}
\label{SI:generalization_scaling}

\textcolor{red}{Here we present results of additional cases in this section to demonstrate the generalizability and scalability of HINTS as discussed and summarized in the main text.}

\textcolor{red}{\textbf{1D Poisson equation with out-of-distribution testing input}. We examine the generalization capability of HINTS in terms of the correlation lengths of the input functions $k(x)$ and $f(x)$, see Fig.~\ref{fig:results_1P_gen_corlen}(D). Different correlation lengths for training data (distinguished by colors of curves) and testing data (horizontal axis) and different volumes of the training dataset (solid/dashed curves) are considered. It is found that the trained HINTS can well generalize to testing functions with correlation lengths higher than the training correlation length; for lower correlation lengths, the convergence rate is slightly worse, but still significantly faster than other methods presented in Fig.~\ref{fig:results_1P}. Such a behavior may be explained by the eigenmode analysis of HINTS presented in Section~\ref{SI:eigenAnalysis}. DeepONet tends to prioritize learning the solutions of low-frequency modes, the spectrum of which is determined by the correlation length of the training data. Test data with larger correlation lengths typically exhibit a narrower span of frequencies, for which the solution has been well learned in the training process. As a consequence, testing with higher correlation lengths results in a similar performance, compared to testing with the same correlation length as training data. In contrast, test data with smaller correlation lengths contains a broader range of high-frequency components, which were not sufficiently represented by the training data. This underrepresentation inevitably results in a degradation in the test performance.}

\textcolor{red}{\textbf{3D Helmholtz equation, large-scale system, tested with a modified computational domain with an internal scatterer}. See Fig.~\ref{fig:results_3H_GS}(B-C). Using this example, we demonstrate (1) the scalability of HINTS to large systems ($N\sim\mathcal{O}(10^6)$, (2) the generalization capability of HINTS in terms of computational domains. The experiment we conduct uses the same trained DeepONet for the 3D Helmholtz problem in with $\approx2\cdot10^3$ parameters, and apply it to a $10^6$ parameters system. In addition, we test the HINTS with an altered problem with the addition of a small cube scatterer inside the 3D domain, without re-training the DeepONet (i.e., the DeepONet is trained on a domain without a scatterer). Technically, the values of $k(\xvec)$ and $f(\xvec)$ for $\xvec$ inside the scatterer are simply padded with zeros. The results show that the HINTS is able to converge, despite the change of the computational domain.}

\textcolor{red}{\textbf{3D Poisson equation, \oldtext{large-scale system} \revv{HINTS vs. direct method}}. To show that the HINTS can solver large scale problems, but can also solve them faster than a direct method, we experiment with the 3D Poisson example with $\mathcal{O}(10^5)-\mathcal{O}(10^6)$ elements. The results for a series of experiments with systems of sizes $6.4\times10^4$, $2.16\times10^5$, and $5.12\times10^5$ elements are shown in Fig.~\ref{fig:results_1P_gen_corlen}(B), where the horizontal axis is the size of the system, and the vertical axis is the time consumption of the cases. We find that HINTS is faster than the direct method for the large system. Furthermore, the logarithmic plot shows that HINTS displays a less steep slope, indicating a better scaling behavior of HINTS compared to the direct method. 
For small-scale systems presented in the main text (1D Poisson), see Fig.~\ref{fig:results_1P_gen_corlen}(A) for a comparison of computational time.}

\textbf{Computational cost of the current implementation of HINTS preconditoner}
The large scale experiment reported in Section~\ref{sec:krylov} is performed using PETSc (version 3.20.4) and Pytorch (version 2.2.0) in serial, using an Xeon E5-2690 v3 processor (64 GB) for numerical linear algebra. 
DeepONet-based components are executed using an NVIDIA Tesla P100 GPU (16 GB).

The run time of our HINTS-based preconditioner implemented in PETSc heavily depends on three factors:
i) a choice of interpolation operator used to map residual to the dimension of the branch network responsible for encoding the right-hand side,
ii) DeepONet architecture, and
iii) the type of GPU used to perform the DeepONet inference.
In our current implementation, most computational time is spent on interpolating the residual between the non-structured finite-element mesh at the finest level of the multigrid hierarchy and the structured mesh used for discretizing branch input functions.
Consequently, to further reduce the runtime of our HINTS solver, we would have to replace currently employed interpolation routines for non-nested meshes, provided by the Firedrake finite element package~\cite{rathgeber2016firedrake}, with more efficient implementations, such as those described in~\cite{farrell2009conservative, krause2016parallel}.
Moreover, the runtime of HINTS-based solvers/preconditioners can be further decreased by considering lighter-weight DeepONet architectures, e.g., latent DeepONet~\cite{kontolati2023learning}, or DITTO~\cite{ovadia2023ditto}, and by utilizing the latest GPU devices.

\textcolor{red}{\textbf{2D Poisson equation, tested with a modified computational domain with an internal cutout and modified boundary conditions.} In another of our recent studies~\cite{kahana2023geometry}, we find that a HINTS trained with PDEs defined in a domain with Dirichlet boundary condition may be used to solve PDEs defined in another domain with mixed Dirichlet and Neumann boundary conditions. This study demonstrates the generalization capability of HINTS in terms of both the computational domain and the boundary conditions. 
In Section~\ref{SI:generalizability}, we provide a preliminary analysis exploring the underlying factor contributing to the generalization capability across computational domains.}

\subsubsection{HINTS for Boundary Conditions Other than Dirichlet}
\label{sec:BC_no_dirichlet}
While we focused on zero Dirichlet boundary conditions in all the examples shown in the main text, we point out that the utility of HINTS is beyond Dirichlet. First, as shown in another of our recent studies~\cite{kahana2023geometry}, HINTS is capable of handling Neumann boundary conditions. In addition, here we show that the HINTS can be used with Sommerfeld radiation boundary conditions, using the case of 2D Helmholtz problem. See Fig.~\ref{fig:results_2H_sommerfeld}(B) for the results. Using Jacobi alone, this problem may not be solved.

\subsubsection{Replacing DeepONet by K-Nearest Neighbors}
\label{sec:k_nearest_neigh}

\textcolor{red}{
To demonstrate the advantages of the DeepONet in HINTS as the operator approximator, we implemented a modified ``HINTS" with the k-nearest neighbor algorithm with Gaussian kernel as the operator approximator. We found that the convergence rate of the alternative solver is $\mu=2.9\times 10^{-2}$, which is significantly smaller than the results for HINTS. We note that the superior predictive capability of DeepONet for (nonlinear) operators expedites the convergence of HINTS. An additional drawback for the alternative method is that the complete training data must be retained for subsequent inference tasks, which poses potential challenges in terms of scalability and resource efficiency.
}

\subsubsection{Preliminary Analysis of Geometry Generalizability.}
\label{SI:generalizability}
\begin{figure}[H]
  \centering
  \includegraphics[width=\textwidth]{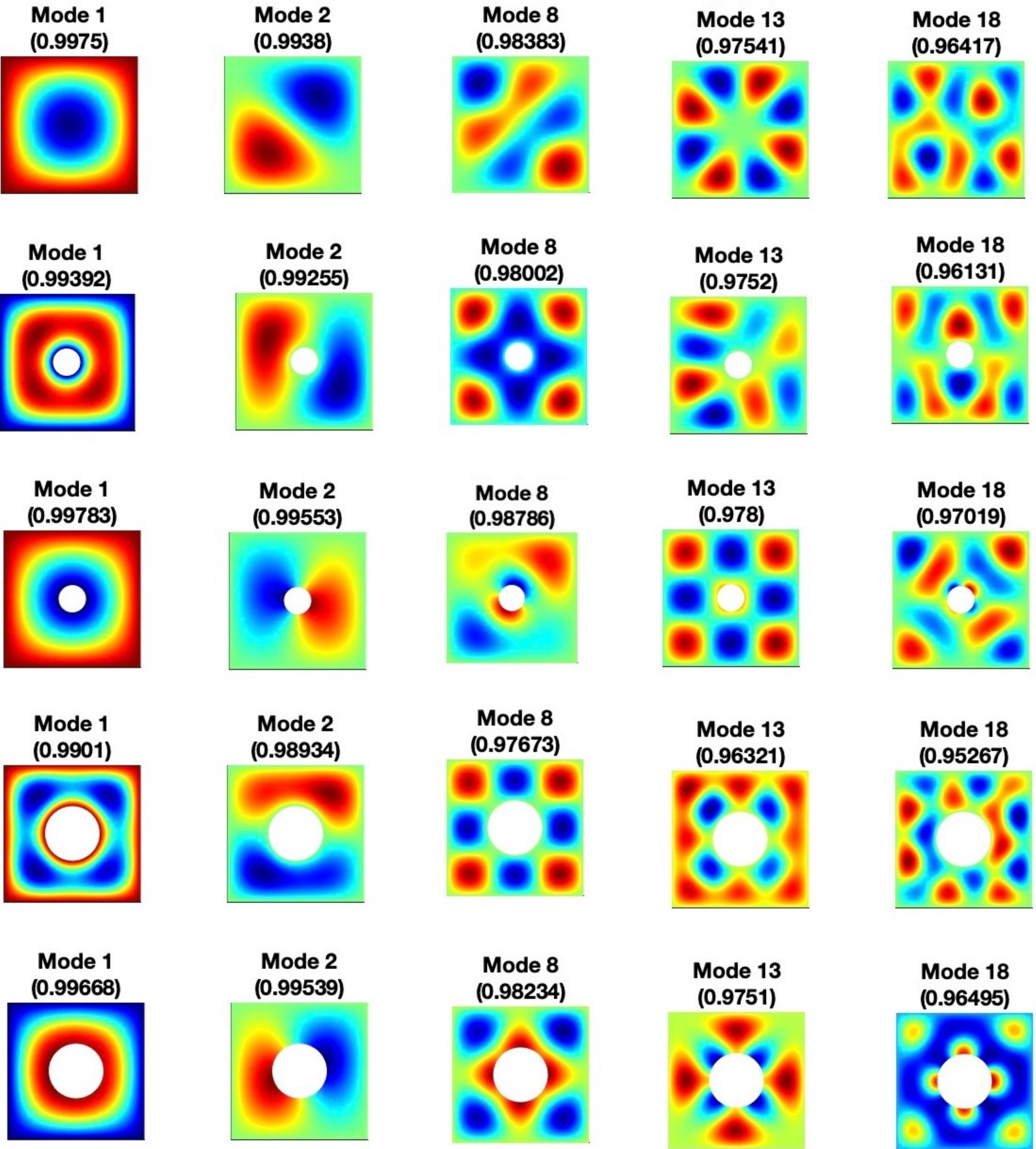}
\caption{\textbf{Eigenmode Analysis of 2D Poisson Equation Solved with HINTS-Jacobi.} 
The number in the parenthesis refers to the eigenvalue corresponding to the eigenmode.
The domain is a unit square and $k(x,y)=1$. 
First row: No cutoff of domain is considered. 
Second/Third row: A circular cutout with radius $0.1$, where Dirichlet/Neumann boundary conditions are applied. 
Fourth/Fifth row:  A circular cutout with radius $0.2$, where Dirichlet/Neumann boundary conditions are applied. }
  \label{fig:Fig_2P_eig_new}
\end{figure}

\textcolor{red}{Here we provide a preliminary analysis on the mechanism of the geometry generalizability of HINTS, which has been presented in Fig.~\ref{fig:results_3H_GS}(B,C) and studied in~\cite{kahana2023geometry}. 
We first consider the 2D Poisson equation defined in a unit square with $k(x,y)=1$. The eigenmodes of the solution are presented in Fig.~\ref{fig:Fig_2P_eig_new}. Then, we include a centered circular cutout into the domain with radius $0.1$ and $0.2$, where Dirichlet or Neumann boundary conditions are enforced on the boundary of the circle (radius $0.1$ or $0.2$)
We find that the low-frequency modes are similar between the cases with and without the introduction of circular cutouts. Because the DeepONet is aimed at predicting the low-frequency modes, such a similarity makes it possible for the HINTS trained in the domain without the cutout to function in the domain with the cutout. This explains why the HINTS is generalizable in terms of the computational domain.
}

\section*{Data Availability}
The datasets used in the manuscript are generated using the provided code; see Section S2.1 for parametric details.  
The dataset for generating the large-scale results (Helmholtz - annular cylinder example) has been uploaded to the open-source Zenodo repository and can be freely accessed at  \\
\href{https://doi.org/10.5281/zenodo.10904349}{https://doi.org/10.5281/zenodo.10904349.}

\section*{Code Availability}
The code used for generating all numerical experiments is publicly available; see the GitHub repository~\cite{hints_precond_code}.

\section*{Acknowledgements}
This work is supported by the DOE PhILMs project (No. de-sc0019453), the MURI-AFOSR FA9550-20-1-0358 projects. 
G.~E.~K. is supported by the ONR Vannevar Bush Faculty Fellowship (N00014-22-1-2795).
A.~Ko.~acknowledges support of the Swiss National Science Foundation (SNF) through the project “Multilevel training of DeepONets – multiscale and multiphysics applications” (206745).

\section*{Author Contributions Statement}
G.~E.~K., J.~P., and E.~T. designed the study and supervised the project. 
E.~Z.~and A.~Ka. developed the method, implemented the computer code, and performed computations. 
A.~Ko developed the PETSc HINTS code, extended the HINTS methodology to preconditioning settings, and designed and performed large-scale experiments.
All authors analyzed the results and contributed to the writing and revising of the manuscript.

\section*{Competing Interests Statement}
G.~E.~K. holds a small equity in Analytica, a private startup company developing AI software products for engineering. 
He provides technical advice on the direction of machine learning to Analytica. 
Analytica has licensed IP from his research related to Physics-Informed Neural Networks.
A.~Ka. and G.~E.~K. are the founders of Phinyx AI, a private startup company developing AI software products for engineering. 
The remaining authors declare no competing interests.

\end{document}